\documentclass {amsart}

\usepackage {hyperref}

\usepackage{amssymb}
\usepackage{amsmath}
\usepackage{amsthm}

\usepackage{enumerate}
\usepackage{enumitem}

\usepackage {comment}

\makeatletter
\@namedef{subjclassname@2010}{%
  \textup{2010} Mathematics Subject Classification}
\makeatother

\newtheorem {theorem} {Theorem} [section]

\newtheorem {lemma} [theorem] {Lemma}
\newtheorem {proposition} [theorem] {Proposition}
\newtheorem {definition} [theorem] {Definition}
\newtheorem {corollary} [theorem] {Corollary}

\newcommand {\apclass} [1] {\ensuremath{\mathrm A_{#1}}}

\newcommand {\lclass} [2] {\ensuremath{\mathrm L_{#1} \left( #2 \right) }}
\newcommand {\lsclass} [1] {\ensuremath{\mathit l^{#1} }}

\newcommand {\hclass} [2] {\ensuremath{\mathrm H_{#1} \left( #2 \right) }}
\newcommand {\lclassg} [1] {\ensuremath{\mathrm L_{#1}}}

\newcommand {\hclassg} [1] {\ensuremath{\mathrm H_{#1}}}

\newcommand {\BMO} {\ensuremath {\mathrm {BMO}}}

\newcommand {\AK} {\ensuremath {\mathrm {AK}}}

\newcommand {\nplus} {\ensuremath {\mathrm {N}^+}}

\DeclareMathOperator* {\supp} {supp}

\newcommand {\weightu} {\ensuremath {\mathit u}}
\newcommand {\weightv} {\ensuremath {\mathit v}}
\newcommand {\weightw} {\ensuremath {\mathit w}}

\newenvironment{titleblock}{}{}
\newenvironment{bibliographyblock}{}{}

\begin{document}
\begin {titleblock}

\baselineskip=17pt


\title {Real interpolation of Hardy-type spaces and BMO-regularity}

\author [D.~V.~Rutsky] {Dmitry V. Rutsky}
\address {
St. Petersburg Department of V.A.Steklov Institute of Mathematics,
Fontanka 27, 191023 St. Petersburg, Russia
}
\email {rutsky@pdmi.ras.ru}

\date{}

\begin{abstract}
Let~$(X, Y)$ be a couple of quasi-Banach lattices of measurable functions
on~$\mathbb T \times \Omega$ satisfying some additional assumptions. The $K$-closedness of a couple of Hardy-type spaces~$(X_A, Y_A)$
in~$(X, Y)$ and the stability of the real interpolation~$(X_A, Y_A)_{\theta, p} = (X_A +
Y_A) \cap (X, Y)_{\theta, p}$ are shown to be equivalent to each other and to
the $\BMO$-regularity of the associated lattices~$\left(\lclassg {1}, \left(X^r\right)' Y^r\right)_{\delta, q}$.
The inclusion~$\left(X^{1 - \theta} Y^\theta\right)_A \subset
  \left(X_A, Y_A\right)_{\theta, \infty}$ is also characterized in these therms.
New examples of couples $(X_A,
Y_A)$ with this stability are given, proving that this property is strictly
weaker than the $\BMO$-regularity of~$(X, Y)$.
\end{abstract}



\thanks{This research was supported by the Russian Science Foundation
(grant No. 18-11-00053).}


\maketitle

\end{titleblock}

\section {Introduction}
\label {introductionS}

This work is mostly concerned with the stability of the real interpolation for
Hardy-type spaces; see \cite {kisliakov1999} for a comprehensive introduction to
the subject with the appropriate references.  In this brief informal
introduction merely an attempt is made to describe a somewhat complicated
subject, the stability of interpolation of Hardy-type spaces and related properties, before going into
detail.

To
give an illustrative example, for classical Hardy spaces the stability is
the formula $$
\left(\hclassg {p}, \hclassg {q}\right)_{\theta, r} = \left[\hclassg {p} +
\hclassg {q}\right] \cap \left(\lclassg {p}, \lclassg {q}\right)_{\theta, r} =
\hclassg {r} $$
with $0 < \theta < 1$, $0 < p < q \leqslant \infty$, $\frac 1 r = \frac {1 -
\theta} p + \frac \theta q$.
For weighted Hardy spaces the corresponding formula
$$
\left(\hclass {p} {\weightu}, \hclass {q} {\weightv}\right)_{\theta, r} =
\left[\hclass {p} {\weightu} + \hclass {q} {\weightv}\right] \cap \left(\lclass
{p} {\weightu}, \lclass {q} {\weightv}\right)_{\theta, r} = \hclass {r} {\weightu^{1- \theta}
\weightv^\theta}
$$
with  $0 < \theta < 1$, $0 < p, q \leqslant \infty$, $\frac 1 r = \frac
{1 - \theta} p + \frac \theta q$
is completely characterized by the condition~$\log \frac \weightu \weightv
\in \BMO$.  Here~$\weightu$ and~$\weightv$ are some suitable weights, and the
weighted spaces are defined by~$X (\weightw) = \weightw X$ with the
corresponding definition for weighted Hardy spaces.

The real and complex interpolation of Hardy spaces were developed in the early
80s, and it was already rather well understood by the mid-90s, including the
weighted and vector-valued cases.
Moreover, for couples~$\left(\hclass {p} {\weightu}, \hclass
{q} {\weightv}\right)$ with~$1 \leqslant p, q \leqslant \infty$ the condition~$\log
\frac \weightu \weightv \in \BMO$ completely characterizes the existence of a partial
retraction from~$\left(\lclass {p} {\weightu}, \lclass
{q} {\weightv}\right)$ and thus also the stability with respect to all
interpolation functors (see~\cite {kisliakov1999}, \cite {kisliakovxu2000};
we note that the existence of a partial retraction and the stability with
respect to general interpolation functors is still unclear for~$p < 1$).

However, interesting questions go beyond the stability for couples of weighted
Hardy spaces.  For spaces~$X$ of functions on the unit circle~$\mathbb T$
the Hardy-type spaces~$X_A$ consist of the functions from suitable spaces of
analytic functions on the unit disk~$\mathbb D$
such that their boundary values are in~$X$.
We work with rather general normed lattices of
measurable functions~$X$, also called the ideal spaces, because they are
characterized by the inclusion~$\lclassg {\infty} X \subset X$.  Some of the
other names for these spaces in wide use are the K\"othe function spaces and
(simply) Banach functional spaces.

At the first glance it may seem as though not much can be done in such a general
setting.  However,
for the complex interpolation of Hardy-type spaces a so-called $\BMO$-regularity
property that generalizes the condition~$\log \frac \weightu
\weightv \in \BMO$ above turned out to be necessary as well as sufficient under rather general
assumptions; see~\cite {kalton1994}.
Intuitively, the stability of interpolation imposes a restriction on the
spaces to be ``smooth'' enough to allow the appropriate decompositions in spaces of
analytic functions.
The $\BMO$-regularity property for one lattice characterizes the smoothness of
its unit ball in a rather peculiar sense, describing the situation when any
given function can be dominated by a weight~$w$ satisfying~$\log w \in \BMO$ with the control on the norm of the
majorant and the~$\BMO$ constant.
These impressive results motivated a substantial amount of research trying, on
one hand, to better understand the $\BMO$-regularity property, and, on the other
hand, to see if this property also characterizes the stability for the real
interpolation and some other interesting related phenomena.

Indeed, for arbitrary couples~$(X_A, Y_A)$ of Hardy-type spaces it was quickly
realized that the $\BMO$-regularity property is sufficient for the stability of
the real interpolation, and even for the $K$-closedness of the corresponding couple.
The $K$-closedness property of~$(X_A, Y_A)$ in~$(X, Y)$ means that arbitrary measurable
decompositions of functions from~$X_A + Y_A$ into a sum of functions from~$X$ and~$Y$ can be
made analytic
with the appropriate control on the norm of the individual parts.
Following~\cite {kisliakov2003}, we call this property the $\AK$-stability of
the couple~$(X, Y)$.

The $K$-closedness property, along with the stability of
interpolation, found many applications in analysis.
We only mention briefly a couple of them here that are most familiar to the
author and illustrate the power of these methods.

In~\cite {kisliakov1989}, a simple proof was found of a rather famous
result~\cite {bourgain1984} that the Grothen\-dieck theorem about $2$-summing
operators holds true for the disc algebra~$C_A$.  See also~\cite
{kisliakov1991en} and~\cite [Chapter~16] {hotgobs1}. This result is derived,
essentially, from the stability of interpolation for couples
of weighted Hardy spaces~$\left(\hclass {2} {\weightw^{-\frac 1 2}}, \hclassg {\infty} \right)$
with certain weights~$\weightw$.
Some related results use the stability of interpolation for
couples of vector-valued Hardy spaces~$\hclass {\infty} {\lsclass {s}}$ and
other subtle interpolation properties.

The results mentioned above appear to be some of the earliest applications of
the stability of interpolation for Hardy spaces, and the constructions involved in
the stability are quite elementary.
These very applications motivated significantly the development
of the theory of interpolation for Hardy-type spaces.  On the other hand, we
mention a recent result~\cite {kisliakov2015pen}, where an old problem about the
vector-valued corona theorem with data in~$\hclass {\infty} {\lsclass {1}}$ was solved with the help of the $\AK$-stability of a
couple of weighted vector-valued spaces~$\left(\lclass {\infty} {\lsclass
{\infty}} (v), \lclass {\infty} {\lsclass {1}} \left(v^{-1}\right)\right)$.
The weight~$v$ arises as a $\BMO$-majorant of a function in~$\lclass
{\infty} {\lsclass {2}}$.  Both the $\AK$-stability of this couple and the
$\BMO$-regularity of~$\lclass {\infty} {\lsclass {2}}$ are rather nontrivial
properties that were only established due to gradual and systematic development
of the theory, mostly during the course of the 90s culminating with~\cite
{kisliakov2002en}, and, at least at present, unlike the former examples they are
{\em not} easily verified by elementary constructions.

A natural question is whether $\BMO$-regularity is also necessary
as well as sufficient for the stability of the real interpolation, and for the
$\AK$-stability.  This problem remained open, at least beyond some special
cases, and the theory of the real interpolation of Hardy-type spaces does not
seem to be complete without a more or less definitive answer to this
question, which remained elusive for the 25 years since the
$\BMO$-regularity property was first introduced in~\cite {kalton1994}.
In~\cite {rutsky2013ben} the author claimed to have proven the equivalence
under some restrictions.  Unfortunately, it was recently discovered that these
results are flawed, and the specific mistake is the formula~\cite [Proposition~18] {rutsky2013ben} $
((X, Y)_{\alpha, p}, Z)_{\beta, p} = (X, (Y, Z)_{\gamma, p})_{\delta, p}
$,
which is false, e.~g., with~$X = \lclassg {2, \infty}$, $Y = \lclassg {1}$, $Z =
\lclassg {\infty}$ and~$\gamma = \frac 1 2$.
It {\em was} established correctly, however, that under some assumptions the
necessary condition for a couple~$(X, Y)$ to be $\AK$-stable is the $\BMO$-regularity of a real
interpolation space~$(\lclassg {1}, X' Y)_{\theta, p}$, and the flawed part is
the derivation of the $\BMO$-regularity of~$X' Y$ from this property.  In the
present work we will see that, surprisingly, these two properties are in fact
{\em not} equivalent, but the former property {\em is} equivalent to
the $\AK$-stability of the couple $(X, Y)$.

The main goal of the present work is to describe comprehensively
the relationship between $\AK$-stability, $\BMO$-regularity and the stability
of the real interpolation of Hardy-type spaces.
We will show that under some
standard assumptions the stability
\begin {equation*}
\label {e:xayast}
(X_A, Y_A)_{\theta, r} = \left[X_A + Y_A\right] \cap (X, Y)_{\theta, r},  
\end {equation*}
the $\AK$-stability of~$(X, Y)$
and the inclusion $\left(X^{1 - \theta} Y^\theta\right)_A \subset \left(X_A,
Y_A\right)_{\theta, \infty}$ are completely characterized by a weaker
$\BMO$-regularity property.
We will call it the {\em weak-type BMO-regularity} property: the
$\BMO$-regularity of a real interpolation space~$\left(\lclassg {1}, (X^r)' Y^r\right)_{\theta,
s}$ with some~$r > 0$, $0 < \theta < 1$ and~$0 < s \leqslant \infty$.
We will see that this property is also equivalent to the $\BMO$-regularity of the
couple $((X, Y)_{\alpha, p}, (X, Y)_{\beta, q})$ with some~$\alpha \neq \beta$
(equivalently, with all~$0 < \alpha < \beta < 1$).

This result leads to previously unknown examples of couples of $\AK$-stable
lattices that are not $\BMO$-regular.  It suffices to consider lattices that coincide
with the Lorentz spaces~$\lclassg {p, q_j}$ with different values of~$q_j$ when
restricted on different sets. In turn, this indicates that even though the
$\BMO$-regularity and weak-type $\BMO$-regularity (equivalent to the $\AK$-stability by our results) are very similar, there is a crucial difference between the two:
the former is stable under multiplication with $\BMO$-regular couples, whereas
the latter is not.

A surprising new observation leading to these results is that a
couple~$(\lclassg {\infty}, Z^\alpha)$ is $\AK$-stable merely
if a related couple~$\left(\lclassg {\infty}, (\lclassg {\infty},
Z)_{\beta, \infty}\right)$ is $\AK$-stable with some~$0 < \alpha < \beta < 1$.
This suggests the aforementioned examples, and it already allows us to
characterize the $\AK$-stability in terms of weak-type $\BMO$-regularity under
the assumption that~$X' Y$ is a Banach lattice.

This observation suggests that, unlike $\BMO$-regularity, the
$\AK$-stability property may be extended to wider couples on a scale.
And indeed, this turns out to be the case: at least under some rather general
assumptions, if~$(E, F)$ is merely an $\AK$-stable couple of lattices of types
$\mathcal C_{\theta_j} (X, Y)$ with some~$0 < \theta_0 < \theta_1 < 1$ then $(X,
Y)$ is also an $\AK$-stable couple.  Thus, the converse is true to a rather
well-known result that couples of real interpolation spaces constructed from an
$\AK$-stable couple are also $\AK$-stable.
This result is rather involved, and unlike all the important results of the
theory up to this point that take advantage of a fixed point
theorem, the Fan--Kakutani fixed point theorem does not suffice for the argument,
and we instead rely on the Powers fixed point theorem for compositions of
acyclic maps.

This indicates that, unlike $\BMO$-regularity,
the $\AK$-stability property of a couple is insensitive to relatively subtle
nuances that do not significantly affect at least some intermediate
couples of spaces of types~$C_{\theta_j}$.
As a consequence, we prove that as soon
as one such couple~$(E, F)$ is $\AK$-stable then all of them are.  In other
words, either all such couples are simultaneously $\AK$-stable, or all of them
are simultaneously not $\AK$-stable. Such a property seems to be rare for the
stability of interpolation of subspaces in general, since it is easy to give examples of one-dimensional subspaces where it fails.

These results depend in a crucial way on a stronger version of $\AK$-stability
that we call the bounded $\AK$-stability, meaning that the respective
analytic decompositions~$H = F + G$ for given measurable decompositions~$H = f +
g \in X_A + Y_A$ can be made by multiplication with some bounded analytic
functions~$U$ and~$1 - U$, i.~e.~$F = U (f + g)$ and~$G = (1 - U) (f + g)$,
which is equivalent to separate control of the norms of~$u = U g$ in~$X$ and~$v
= (1 - U) f$ in~$Y$.
These functions~$u$ and~$v$ thus belong to the intersection~$X \cap Y$, so their
norms can also be usefully estimated in terms of the norms of the spaces
of type~$\mathcal C_\theta (X, Y)$.

It was noticed
in~\cite {rutsky2010en} that the bounded $\AK$-stability property allows us to
improve the convexity of $\AK$-stable lattices.
We further develop these ideas, and prove a
rather general result showing how even countable analytic decompositions, such
as those that arise in the real interpolation spaces of couples of Hardy-type
spaces, can often be made bounded in a similar sense.  This allows us to improve
the convexity of the couples that are stable with respect to the real interpolation, and gain the
convexity that is required to apply the results~\cite {kalton1994}
on the stability of the complex interpolation.  We mention that there
is a different approach to~\cite {kalton1994} to be published elsewhere that
altogether avoids the nontrivial convexity assumptions.  However, it
still takes advantage of the general result about bounded analytic
decompositions in the present work.

The plan of the paper is as follows.  In the next Section~\ref {s:somr} we
give all the necessary formal definitions and state the main results. 
At the end of it we present examples of couples of lattices that are $\AK$-stable but not
$\BMO$-regular.
In Section~\ref {s:spotwtbr} we establish some useful properties of the
weak-type $\BMO$-regularity that mostly follow the properties of the usual
$\BMO$-regularity.
Then, in Section~\ref {basS} the properties of the bounded $\AK$-stability are
studied, and we prove the main results for the case of nondiscrete additional
variable.  These results are much weaker compared to what we can do in the
discrete case, but they are also rather short and simple, do not use fixed
point arguments, and give a good idea of what is going on in the discrete case.

Most results after that require the second variable to be discrete.
In Section~\ref {akbakeqs} we briefly discuss the phenomenon of bounded
$\nplus$-stability and bounded analytic decompositions in general, which also
covers the notion of bounded $\AK$-stability, and state a rather general result
which will be proven in Section~\ref {s:jxsp} with the help of a fixed point theorem. 
In order to do this, in Section~\ref {s:tuccs} the topology of pointwise
convergence on compact sets is introduced for Hardy-type spaces~$X_A$, and we prove that the Fatou property of the
lattice~$X$ implies that the unit ball of~$X_A$ is compact with respect to this
topology.  As an application, we show that the so-called strong $\AK$-stability
is equivalent to the usual one for quasi-Banach lattices.

The remaining thee sections contain the proof of the main result.  In
Section~\ref {s:socxy} the $\AK$-stability of a couple~$(X, Y)$ is derived from
the inclusion $\left(X^{1 - \theta} Y^\theta\right)_A \subset \left(X_A,
Y_A\right)_{\theta, \infty}$.  This also provides a nice, short and
essentially self-contained proof that the stability for the real interpolation
functor~$(\cdot, \cdot)_{\theta, \infty}$ is equivalent to the $\AK$-stability
for a Banach couple with the Fatou property.
Moreover, the much simpler Fan--Kakutani theorem can be used in this
argument instead of the Powers theorem with little additional effort. 
Section~\ref {s:trtbakst} contains the proof of the crucial observation that the
bounded $\AK$-stability of a couple of spaces of type~$\mathcal C_{\theta_j} (X, Y)$ implies the bounded $\AK$-stability of~$(X,
Y)$.  Finally, in Section~\ref {s:goodrealinterp} we prove the main result.

\section {Statement of the main results}

\label {s:somr}

We mostly work with spaces of measu\-rable functions
on the measurable space $\mathbb T \times \Omega$,
where $\mathbb T$ is the unit circle with the Lebesgue measure and $(\Omega,
\mu)$ is some $\sigma$-finite measurable space that represents an additional
variable.  For technical reasons we will often assume~$\Omega$ to be discrete,
which means that~$\Omega$ is at most countable.

A quasi-normed lattice of mea\-surable functions~$X$ on a
$\sigma$-finite measurable space~$\mathcal M$ is a quasi-normed space of
measurable functions~$X$ in which the norm is compatible with the natural order; that is, if~$|f| \leqslant g$ for some function~$g \in X$ then~$f \in X$
and~$\|f\|_X \leqslant \|g\|_X$.
For simplicity we only work with lattices~$X$ such that~$\supp X = \mathcal M$
up to a set of measure~$0$.
For more detail on the normed lattices and their properties see, e.~g., \cite [Chapter~10] {kantorovichold}.

We say that~$X$ has the Fatou property if for any~$f_j \in X$,
$\|f_j\|_X \leqslant 1$ such that~$f_j \to f$ almost everywhere we also have~$f
\in X$ and~$\|f\|_X \leqslant 1$.  For normed lattices~$X$ the Fatou property is
equivalent to the closedness of the unit ball~$B_X$ with respect to the
convergence in measure on sets of finite measure, and it implies that~$X$ is a
Banach lattice.
$X$ is said to have order continuous norm if for any nonincreasing sequence~$f_j
\in X$ converging to~$0$ we also have~$\|f_j\|_X \to 0$.  The order dual~$X'$
can be defined as a lattice with the norm~$\|g\|_{X'} = \sup_{f \in B_X} \int |f
g|$.  The Fatou property is equivalent to~$(X')' = X$, and the order continuity
is equivalent to~$X' = X^*$.  For example, $\lclassg {p}' = \lclassg {p'}$ for
all~$1 \leqslant p \leqslant \infty$.

For quasi-normed lattices~$X$ and~$Y$ of measurable functions on a
$\sigma$-finite measurable space~$\mathcal M$ the pointwise product~$X Y$
is defined by the quasinorm~$\|h\|_{X Y} = \inf_{h = f g} \|f\|_X \|g\|_Y$, and
the power~$X^\delta$, $\delta > 0$ (sometimes called the $\frac 1
\delta$-convexification of~$X$) is defined by the quasinorm~$\|f\|_{X^\delta} =
\left\| |f|^{1 \slash \delta}\right\|^\delta$.
This allows us to define a Calder\'on-Lozanovski\u{\i} product $Z = X^{1 -
\theta} Y^\theta$, $0 < \theta < 1$, which naturally inherits many properties from~$X$
and~$Y$.  If~$X$ and~$Y$ are Banach
lattices with the Fatou property then so is the product~$Z$, the dual can be
computed as $Z' = X'^{1 - \theta} Y'^\theta$, and~$\lclassg {1} = X' X$ by the
Lozanovski\u{\i} factorization theorem (see \cite {lozanovsky1969}).

Let~$\nplus$ be the set of boundary values of the Smirnov class of analytic
functions on the disc (see, e.~g., \cite {privaloven}, \cite {hoffman}).  By~$\nplus \otimes \Omega$ we
understand the set of measurable functions~$f$ on~$\mathbb T \times \Omega$ such
that~$f (\cdot, \omega) \in \nplus$ for almost all~$\omega \in \Omega$.  For a
space~$X$ of measurable functions on~$\mathbb T \times \Omega$ we define the
corresponding Hardy-type space~$X_A = X \cap (\nplus \otimes \Omega)$.
For example, from the Lebesgue spaces~$\lclassg {p}$, $0 < p \leqslant \infty$
we get the usual Hardy spaces $\left[\lclassg {p}\right]_A = \hclassg {p}$, but
this definition also yields the Hardy-Lorentz spaces $\hclassg {p, q}$, the weighted Hardy spaces $\hclass {p} {\weightw}$,
the variable exponent Hardy spaces $\hclassg {p (\cdot)}$, the vector-valued Hardy spaces $\hclass {p} {\lsclass {q}}$
and many others.

\begin {definition}
\label {d:starpr}
Suppose that~$X$ is a quasi-normed lattice of measurable functions on~$\mathbb T
\times \Omega$.  We say that~$X$ satisfies property $(*)$ with constant~$C$ 
if for any~$f \in X$, $f \neq 0$ there exists a majorant~$g \geqslant |f|$
such that~$\|g\|_X \leqslant C \|f\|_X$ and $\log g (\cdot, \omega) \in \lclassg
{1}$ for almost all~$\omega \in \Omega$.
\end{definition}
This property is often assumed to avoid degeneration.  It says, essentially,
that lattice~$X$ has a complete set of outer functions.
If it is satisfied, then by~\cite [Lemma~2.2] {kalton1994}
it is also satisfied with arbitrary constants~$C > 1$.

Let~$r > 0$.  If a $r$-convex quasi-normed lattice~$X$ of measurable functions
on~$\mathbb T \times \Omega$ has the Fatou property and property~$(*)$
then~$X_A$ is a closed subspace of~$X$; see, e.~g., \cite [\S 1] {kisliakov2002en} (where~$X$ is
assumed to be normed, but this is easily generalized to quasi-normed
lattices).

One of the interesting questions of the theory of interpolation spaces is their
stability with respect to the intersection with subspaces.  In the present
work we only consider the stability with respect to the intersection with
spaces of analytic functions, and so we only give the definitions of rather
general phenomena specialized to the case of Hardy-type spaces.
For the interpolation theory see, e.~g., \cite {bergh}.
\begin {definition}
\label {d:nplusst}
We say that a couple~$(X, Y)$ of quasi-normed lattices of measurable functions
on~$\mathbb T \times \Omega$ is $\nplus$-stable with respect to an interpolation
functor~$\mathcal F$ if~$\left[\mathcal F \left((X, Y)\right)\right]_A =
\mathcal F \left((X_A, Y_A)\right)$.
\end{definition}

For the real interpolation functors the $\nplus$-stability is implied by the
following property, which is on its own of considerable interest.
\begin {definition}
\label {akstability}
A quasi-normed couple~$(X, Y)$ of lattices of measurable functions on~$\mathbb T
\times \Omega$ is called $\AK$-stable with constant~$C$ if~$(X_A, Y_A)$ is
$K$-closed in $(X, Y)$ with constant~$C$.  That is, for any~$H \in X_A + Y_A$
and $f \in X$, $g \in Y$ such that $H = f + g$ there exist some $F \in X_A$, $G
\in Y_A$ such that~$H = F + G$ and~$\|F\|_X \leqslant C \|f\|_X$,
$\|G\|_Y \leqslant C \|g\|_Y$.
\end {definition}

The $\BMO$-regularity properties introduced below were found to be closely
related to the above properties.
For the first time they were explicitly introduced and extensively studied,
apparently, in~\cite {kalton1994} in order to characterize the stability of
the complex interpolation for Hardy-type spaces, and then in~\cite
{kisliakov1999} for both the real and the complex interpolation, although they were also somewhat
implicitly used before in a different form (later found to be equivalent to
$\BMO$-regularity) in various stability results such as~\cite {kisliakovxu1994}.
\begin {definition}
\label {bmordef}
A quasi-normed lattice~$X$ of measurable functions on~$\mathbb T \times
\Omega$ is called $\BMO$-regular with constants~$(C, m)$ if for any nonzero~$f
\in X$ there exists a majorant~$u \geqslant |f|$ such that~$\|u\|_X \leqslant m \|f\|_X$ and $\left\|\log u (\cdot,
\omega)\right\|_{\BMO} \leqslant C$ for almost all~$\omega \in \Omega$.
\end {definition}
As a quick example we mention that all rearrangement invariant
lattices that are intermediate spaces for the couple~$(\lclassg {1}, \lclassg
{\infty})$, such as the Lorentz spaces~$\lclassg {p, q}$ with~$1 \leqslant p, q
\leqslant \infty$, are $\BMO$-regular (see, e.~g., \cite [Proposition~2]
{rutsky2011en}).  On the other hand, if~$X$ is a $\BMO$-regular lattice then
the weighted lattice~$X (\weightw)$ is $\BMO$-regular if and only if~$\log
\weightw (\cdot, \omega) \in \BMO$ uniformly in almost all~$\omega \in \Omega$
(see, e.~g, \cite [Proposition~5] {rutsky2011en}).

\begin {definition}
\label {bmorcdef}
A couple~$(X, Y)$ of quasi-normed lattices of measurable functions on~$\mathbb
T \times \Omega$ is said to be $\BMO$-regular with constants~$(C, m)$ if for all
nonzero~$f \in X$ and~$g \in Y$ there exist some majorants~$u \geqslant |f|$ and~$v \geqslant |g|$ such
that~$\|u\|_X \leqslant m \|f\|_X$, $\|v\|_Y \leqslant m \|g\|_Y$ and
$\left\|\log \frac {u (\cdot, \omega)} {v (\cdot, \omega)}\right\|_{\BMO}
\leqslant C$ for almost all $\omega \in \Omega$.
\end {definition}
It is easy to see that if both lattices~$X$ and~$Y$ are
$\BMO$-regular then couple~$(X, Y)$ is also $\BMO$-regular.
If~$X$ and~$Y$ are $r$-convex with some~$r > 0$ then the $\BMO$-regularity
of~$(X, Y)$ is equivalent to the $\BMO$-regularity of $(X^r)' (Y^r)$ for
lattices with the Fatou property (see \cite [Theorem~8] {rutsky2011en}).

It is well known that $\BMO$-regularity of a couple~$(X, Y)$
implies its $\AK$-stability (see, e.~g., \cite [Theorem~3.3]
{kisliakov1999}).  This (up to some detail) follows from the fact that
couples~$(\lclass {\infty} {\weightu}, \lclass {\infty} {\weightv})$ are $\AK$-stable for
the corresponding $\BMO$-majorants $(\weightu, \weightv)$.
The converse was long suspected to be true, i.~e.
that some kind of $\BMO$-regularity is also necessary for $\AK$-stability, and for couples of weighted Lebesgue spaces $\AK$-stability
is indeed equivalent to $\BMO$-regularity (see \cite [Theorem~3.2]
{kisliakov1999}, \cite [Theorem~1] {rutsky2009en} with the original result
obtained in~\cite [Theorem~1.8] {cwikelmccarthywolf1992}).
There are also some couples with additional variable for which it is true
(see~\cite [Theorem~1] {kisliakov2002en} and~\cite [Theorem~2] {rutsky2010en}). 
However, we will show that under some natural assumptions $\AK$-stability and even the stability with respect to the
real interpolation are completely characterized in terms of a weaker property.

\begin {definition}
\label {wtbmor}
Suppose that~$(X, Y)$ is a couple of quasi-Banach latticers of measurable
functions on~$\mathbb T \times \Omega$
such that~$X$ is $r$-convex with some $r > 0$.
We say that~$(X, Y)$ is weak-type $\BMO$-regular if $\left(\lclassg {1},
(X^r)' Y^r\right)_{\theta, p}$ is $\BMO$-regular with some~$0 < \theta < 1$ and
$0 < p \leqslant \infty$.
\end{definition}
This definition is meant to be understood in the sense that the
$\BMO$-regularity is present for small enough values of~$\theta$.
By Proposition~\ref
{wtbmoruniv} below,
Definition~\ref {wtbmor} does
not essentially depend on~$p$ and~$r$.

By~\cite [Theorem~8] {rutsky2011en} mentioned above and Proposition~\ref
{rintbmor} below, a $\BMO$-regular couple of quasi-Banach
lattices is also weak-type $\BMO$-regular.
The converse is false in general; see examples at the end of this section.
One crucial difference to note between the $\BMO$-regularity and weak-type
$\BMO$-regularity is that, unlike the former, the latter is not stable
under the multiplication of couples in its various specific forms, i.~e.
if~$(X, Y)$ and~$(E, F)$ are both weak-type $\BMO$-regular then~$(X E, Y
F)$ is not necessarily weak-type $\BMO$-regular, not even for couples~$(E, F) =
\left(\lclass {p} {\lsclass {p}}, \lclass {q} {\lsclass {q}}\right)$ with~$p
\neq q$.
Otherwise the main result coupled with~\cite [Theorem~2] {rutsky2010en}
would have given us the equivalence of the weak-type $\BMO$-regularity to the
$\BMO$-regularity.  Without the additional variable this multiplication also
fails for $\BMO$-regular couples of weighted Lebesgue spaces~$(E, F) =
\left(\lclass {p} {u}, \lclass {q} {v}\right)$ (see~\cite [Proposition~21]
{rutsky2010en}).
It still is not clear, however, whether the weak-type $\BMO$-regularity is
stable under multiplication by at least a couple of unweighted Lebesgue
spaces~$(E, F) = \left(\lclassg {p}, \lclassg {q}\right)$ without the additional variable.

On the other hand, the distinction between these properties does not appear to
be big.  The equivalence of conditions~\ref {t:c:isbmor} and~\ref
{t:c:xybmor} of Theorem~\ref {goodrealinterp} below shows that under its
assumptions the $\BMO$-regularity is equivalent to the weak-type
$\BMO$-regularity for couples obtained by the real interpolation from a single
couple.
For example, these properties coincide for couples of weighted Lebesgue spaces. 
Since both of these conditions are invariant under raising the lattices to
any positive power, the convexity assumptions in this equivalence may be further
relaxed away to the assumption that both lattices are $r$-convex with some~$r > 0$.

There are other interesting natural spaces to be investigated where
one might suspect the equivalence of the weak-type $\BMO$-regularity to the usual one, such as couples of weighted
vector-valued Lebesgue spaces~$\lclass {p} {\lsclass {q}} (\weightw)$, weighted
Orlicz spaces and variable exponent Lebesgue spaces~$\lclassg {p (\cdot)}$ to
name a few.

Now we are ready to state the main result of the present work.  It establishes
that under some standard assumptions the weak-type $\BMO$-regularity completely
characterizes various properties related
to the stability of the real interpolation.  Moreover, it also shows that these
properties are closely related to one another, and both $\AK$-stability and
weak-type $\BMO$-regularity are invariant under the transition to other
couples on a single interpolation scale.
\begin {theorem}
\label {goodrealinterp}
Suppose that~$(X, Y)$ is a couple of quasi-normed $r$-convex lattices of
measurable functions on~$\mathbb T \times \Omega$ with a discrete space~$\Omega$
and some~$r > 0$ satisfying the Fatou property and property~$(*)$
such that~$X^{1 - \theta_j} Y^{\theta_j}$ are Banach lattices with some $0 <
\theta_0 < \theta < \theta_1 < 1$.
The following conditions are equivalent.
\begin {enumerate} [label=\text {\upshape(\roman*)}, leftmargin=*, widest=iiii]
  \item \label {t:c:npluss} $(X, Y)$ is $\nplus$-stable with respect to~$(\cdot,
  \cdot)_{\theta, s}$ for some (equivalently, for all) $1 \leqslant s \leqslant
  \infty$.
  \item \label {t:c:akst} $(X, Y)$ is $\AK$-stable.
  \item \label {t:c:efakst} $(E, F)$ is $\AK$-stable for some (equivalently, for
  all) quasi-normed $r$-convex lattices~$E$ and~$F$ of measurable functions
  on~$\mathbb T \times \Omega$ satisfying the Fatou property and property~$(*)$
  such that~$E$ is of type~$\mathcal C_{\alpha} (X, Y)$ and~$F$ is of
  type~$\mathcal C_\beta (X, Y)$ with some~$0 \leqslant \alpha < \theta <
  \beta \leqslant 1$.
  \item \label {t:c:isakst} $\left( (X, Y)_{\alpha, p}, (X, Y)_{\beta, q}
  \right)$ is $\AK$-stable with some~$0 < \alpha < \theta < \beta < 1$
  (equivalently, with all~$0 < \alpha < \beta < 1$) and~$0 < p, q \leqslant
  \infty$.
  \item \label {t:c:isbmor} $\left( (X, Y)_{\alpha, p}, (X, Y)_{\beta, q}
  \right)$ is $\BMO$-regular with some (equivalently, with all)~$0 < \alpha < \beta < 1$ and~$0 < p, q
  \leqslant \infty$.
  \item \label {t:c:xybmor} $\left( (X, Y)_{\alpha, p}, (X, Y)_{\beta, q}
  \right)$ is weak-type $\BMO$-regular with some (equivalently, with all)~$0 < \alpha <
  \beta < 1$ and~$0 < p, q \leqslant \infty$.
  \item \label {t:c:efbmor} $(E, F)$ is weak-type $\BMO$-regular for some
  (equivalently, for all) lattices~$E$ and~$F$ defined in condition~\ref
  {t:c:efakst}
  \item \label {t:c:wtbmor} $(X, Y)$ is weak-type $\BMO$-regular.
  \item \label {t:c:xyincl} $\left(X^{1 - \theta} Y^\theta\right)_A \subset
  \left(X_A, Y_A\right)_{\theta, \infty}$.
\end{enumerate}
\end {theorem}
In particular, a couple~$(X, Y)$ satisfying the assumptions of Theorem~\ref
{goodrealinterp} may be $\nplus$-stable with respect to the real
interpolation~$(\cdot, \cdot)_{\theta, q}$ but not with respect
to the complex interpolation~$(\cdot, \cdot)_\theta$, since the latter stability
is equivalent to the $\BMO$-regularity of the couple~$(X, Y)$ (at least under
some mild convextiy assumptions) by~\cite [Theorem~5.12] {kalton1994}.

Condition~\ref {t:c:xyincl} generalizes
the corresponding result for couples of weighted Lebes\-gue spaces (see \cite
[Theorem~1.8] {cwikelmccarthywolf1992}, \cite [Theorem~3.2] {kisliakov1999},
\cite [Theorem~1] {rutsky2009en}).
As a consequence, this shows that if $X^{1 -\theta} Y^\theta$ has order
continuous norm and~$(X, Y)$ is $\nplus$-stable with respect to an interpolation
functor~$\mathcal F \supset (\cdot, \cdot)_{\theta}$ of type~$\mathcal
C_\theta$ then~$(X, Y)$ is weak-type $\BMO$-regular.  This suggests an
interesting question: for which functors~$\mathcal F$ of type~$\mathcal
C_\theta$ the same holds true?  The positive answer for all such~$\mathcal F$ is
equivalent to the statement
that the inclusion~$\left[\left(X, Y\right)_{\theta, 1}\right]_A \subset
\left(X_A, Y_A\right)_{\theta, \infty}$ implies that~$(X, Y)$ is weak-type
$\BMO$-regular.

We note that although the assumptions made in Theorem~\ref {goodrealinterp} are fairly broad,
the generality of these results is still not entirely satisfying.  In
particular, the convexity assumptions on the lattices and the discreteness
assumption on~$\Omega$ arise because they provide crucial geometrical and
topological properties of certain maps used in the proofs, and at present it is
not clear how to get around these restrictions.
A different approach yields the connection between $\AK$-stability and $\BMO$-regularity for
arbitrary~$\Omega$ but with certain restrictions on lattices.
\begin {theorem}
\label {akbmoreq0}
Let~$(X, Y)$ be a couple of quasi-Banach lattices of measurable functions
on~$\mathbb T \times \Omega$ satisfying the Fatou property and property~$(*)$ such that
either~$Y = \lclassg {\infty}$ or~$X$, $Y$ and~$X' Y$ are Banach lattices. 
Then couple~$(X, Y)$ is $\AK$-stable if and only if it is weak-type $\BMO$-regular.
\end{theorem}
The proof of Theorem~\ref {akbmoreq0} is given in Section~\ref
{bakstability} below (see Propositions~\ref {akbmoreq0i} and~\ref {akbmoreq0r}), and on
its own it is fairly uncomplicated.  The individual transitions are obtained
under broader sets of assumptions that all reduce to the case~$Y = \lclassg
{\infty}$.
The ``only if'' part, for example, is valid if~$X = \lclassg {p}$ with any~$1 \leqslant p \leqslant \infty$ and~$Y$ is
a Banach lattice.  The ``if'' part is somewhat less satisfactory, requiring
that~$(X^r)' Y^r$ is a Banach lattice with some~$r > 0$, but this at least
covers all couples of weighted Lebesgue spaces, or, more generally,
couples~$(X, Y)$ such that~$X$ is $q$-concave and~$Y$ is $q$-convex with some~$0
< q \leqslant \infty$.

We now give some examples of lattices~$Y$ such that couple~$(\lclassg {1}, Y)$
is weak-type $\BMO$-regular and hence $\AK$-stable but~$Y$ is not $\BMO$-regular.
Let~$\mu$ be a point mass, i.~e. we do not consider the
additional variable.  For a measurable set~$E \subset \mathbb T$ and
quasi-normed lattices~$Y_0$ and~$Y_1$ of measurable functions on~$\mathbb T$ we
define a composite lattice
$$
Y = \chi_{\mathbb T \setminus E} Y_0 + \chi_E Y_1 =
\left\{\chi_{\mathbb T \setminus E} f + \chi_E g \mid f \in Y_0, g \in
Y_1 \right\}
$$
with a norm~$\|\chi_{\mathbb T \setminus E} f + \chi_E g\|_Y =
\|\chi_{\mathbb T \setminus E} f\|_{Y_0} +\|\chi_{E}
g\|_{Y_1}$.
For simplicity we choose the half-circle~$E = [0, \pi)$ and
the Lorentz spaces~$Y_j =\lclassg
{t, {s_j}}$, $j \in \mathbb \{0, 1\}$ with some~$1 < t < \infty$, $0 < s_0, s_1
\leqslant \infty$.
Then it is easy to see that~$\left( \lclassg {1}, Y \right)_{\theta, p} =
\lclassg {q, p}$ is $\BMO$-regular with~$0 < \theta < 1$ and~$q = \left(\frac {1
- \theta} 1 + \frac \theta t\right)^{-1}$, so~$(\lclassg {1}, Y)$ is weak-type
$\BMO$-regular.  However, $Y$ is not $\BMO$-regular for~$s_0 \neq s_1$.
To see this, suppose that, more generally, $Y_0 \subsetneqq Y_1$ are
some rearrangement invariant spaces with the Fatou property
and~$Y$ is
$\BMO$-regular.
Then by \cite [Theorem~1] {rutsky2011en} the Hardy-Littlewood maximal operator~$M f (e^{i x})
=
\sup_{0 < r \leqslant \pi} \int_{x - r}^{x + r} |f (e^{i s})| ds$, $x \in
\mathbb R$ is bounded in~$Z = \left({\lclassg {1}\strut}^{1 - \alpha}
{Y\strut}^\alpha\right)^\beta$ for some~$0 < \alpha, \beta < 1$.
Observe that~$Z = \chi_{\mathbb T \setminus E} Z_0 + \chi_E Z_1$
with~$Z_j = \left({\lclassg {1}\strut}^{1 - \alpha}
{Y_j\strut}^\alpha\right)^\beta$, $j \in \{0, 1\}$, and~$Z_0 \subsetneqq Z_1$:
otherwise the equality~$Z_0 = Z_1$ would imply that~${\lclassg {1}\strut}^{1 -
\alpha} {Y_0\strut}^\alpha = {\lclassg {1}\strut}^{1 - \alpha} {Y_1\strut}^\alpha$,
$Y_0'^\alpha = \left({\lclassg {1}\strut}^{1 - \alpha}
{Y_0\strut}^\alpha\right)' = \left({\lclassg {1}\strut}^{1 - \alpha}
{Y_1\strut}^\alpha\right)' = Y_1'^\alpha$, $Y_0' = Y_1'$ and~$Y_0 = (Y_0')'
= (Y_1')' = Y_1$.  Let~$f \in Z_1 \setminus Z_0$, and let
$$
g (e^{i x}) =
\chi_{[0, \pi)} (x) \left(t \mapsto \chi_{[0, 2 \pi)} (t) f (e^{i t}) \right)^*
(x), \quad x \in [0, 2 \pi) $$
be the nonincreasing rearrangement of~$f$ restricted to the upper half-circle.
Then~$g \in Z$, $M g \in Z$ and in particular~$\chi_{\mathbb T \setminus E} M g
\in Z_0$. It is easy to see that~$g \in Z_1 \setminus Z_0$ and~$M g (e^{-i x})
\geqslant \frac 1 {4 x} \int_0^x g \geqslant \frac 1 4 g (x)$ for~$0 < x < \frac \pi 2$,
which contradicts~$\chi_{\mathbb T \setminus E} M g \in Z_0$.

\section {Some properties of weak-type BMO-regularity}
\label {s:spotwtbr}

The following formula (also appearing in \cite [Lemma~1] {kisliakov2015pen}
with a short proof) seems to be rather well known; see, e.~g., \cite
[Theorem~3.7] {schep2010}.
\begin {proposition}
\label {xyxm}
Suppose that~$E$ and~$F$ are Banach lattices of measurable functions on the same
$\sigma$-finite measurable space having the Fatou property
such that $E F$ is also a Banach lattice.  Then $E' = (E F)' F$.
\end {proposition}

\begin {proposition} {\cite [Proposition~14] {rutsky2013ben}}
\label {p:riscale}
Let~$X$ and~$Y$ be some quasi-Banach lattices of measurable functions on some
$\sigma$-finite measurable space.
Then
$$
(X, Y)_{\theta, p}^\alpha = \left(X^\alpha, Y^\alpha\right)_{\theta, \frac
p \alpha}
$$ for all $\alpha > 0$, $0 < \theta < 1$ and $0 < p \leqslant \infty$.
\end {proposition}

The following observation is a simple consequence of the well-known fact
that a lattice~$Z$ is~$\BMO$-regular if and only if~$Z^\delta$ is~$\apclass
{2}$-regular, \cite [Proposition~17] {rutsky2013ben} and~Proposition~\ref
{p:riscale}.
\begin {proposition}
\label {rintbmor}
Suppose that quasi-Banach lattices~$X$ and~$Y$ of measurable functions
on~$\mathbb T \times \Omega$ are $\BMO$-regular.
Then the real interpolation space~$(X, Y)_{\theta, q}$ is also a $\BMO$-regular
lattice for all~$0 < \theta < 1$ and~$0 < q \leqslant \infty$.
\end{proposition}

Defition~\ref {wtbmor}	of weak-type $\BMO$-regularity does not
depend on~$r$ and~$p$.
\begin {proposition}
\label {wtbmoruniv}
Suppose that~$(X, Y)$ is a couple of quasi-normed latticers of measurable
functions on~$\mathbb T \times \Omega$
such that~$X$ is $r$-convex with some~$r > 0$.
Then~$(X, Y)$ is weak-type $\BMO$-regular if and only if
$\left(\lclassg {1}, {\lclassg {1}\strut}^{1 - \gamma} {\left[(X^s)'
Y^s\right]\strut}^\gamma\right)_{\theta, q}$ is $\BMO$-regular for some~$0 <
\theta < 1$ (equivalently, for all sufficiently small~$\theta > 0$) and for
some (equivalently, for all)~$0 < s \leqslant r$, $0 < \gamma \leqslant 1$
and~$0 < q \leqslant \infty$.
\end{proposition}
Indeed, by Proposition~\ref {rintbmor} we can take
arbitrary~$p$ in Definition~\ref
{wtbmor} with any smaller~$\theta$, since lattice~$\lclassg {1}$ is
$\BMO$-regular and by the reiteration theorem (see, e.~g., \cite
[Theorem~3.5.3] {bergh}) we have $$
\left(\lclassg {1}, \left(\lclassg {1},
(X^r)' (Y^r)\right)_{\theta, p}\right)_{\eta, q} = \left(\lclassg {1},
(X^r)' (Y^r)\right)_{\eta \theta, q}
$$
for arbitrary~$0 < \eta < 1$ and~$0 < q \leqslant \infty$.
For the independence from~$r$, observe that~$(X^s)' Y^s = ([X^r]^{\frac s r})'
Y^s = {\lclassg {1}\strut}^{1 - \frac s r}{\left[(X^r)' Y^r\right]\strut}^{\frac s
r}$ and~${\lclassg {1}\strut}^{1 - \gamma} {\left[(X^s)'
Y^s\right]\strut}^\gamma = {\lclassg {1}\strut}^{1 - \frac s r
\gamma}{\left[(X^r)' Y^r\right]\strut}^{\frac s r \gamma}$, thus by the
reiteration theorem
$$
\left(\lclassg {1}, {\lclassg {1}\strut}^{1 - \gamma} {\left[(X^s)'
Y^s\right]\strut}^\gamma\right)_{\theta, q} = \left(\lclassg {1}, (X^r)' Y^r\right)_{\theta \frac s
r \gamma, q}.
$$
\begin {corollary}
\label {wtbmounif}
Suppose that~$(X, Y)$ is a couple of quasi-normed latticers of measurable
functions on~$\mathbb T \times \Omega$
such that~$X$ is $r$-convex with some~$r > 0$.  If~$(X, Y)$ is weak-type
$\BMO$-regular then so is~$(X^\delta, Y^\delta)$ for all~$\delta > 0$.
\end{corollary}

The weak-type $\BMO$-regularity has the natural symmetry,
duality and divisibility properties.  In the present work these properties
are only used in the proof of Proposition~\ref {akbmoreq0r} below under the
assumptions unrelated to Theorem~\ref {akbmoreq0}.
\begin {proposition}
\label {wtbmordiv}
Suppose that~$X$, $Y$ and~$Z$ are $r$-convex quasi-normed lattices of
measurable functions on~$\mathbb T \times \Omega$ with some~$r > 0$ satisfying
the Fatou property.  The following conditions are equivalent.
\begin {enumerate} [label=\text {\upshape(\roman*)}, leftmargin=*, widest=iii]
  \item \label {p:c:xywtbmor} $(X, Y)$ is weak-type $\BMO$-regular.
  \item \label {p:c:xzyzwtbmor} $(X Z, Y Z)$ is weak-type $\BMO$-regular.
  \item \label {p:c:yrxrwtbmor} $\left((Y^r)', (X^r)'\right)$ is weak-type
  $\BMO$-regular.
  \item \label {p:c:yxwtbmor} $(Y, X)$  is weak-type $\BMO$-regular.
\end{enumerate}
\end{proposition}
Indeed, by Corollary~\ref {wtbmounif} condition~\ref {p:c:xywtbmor} is
equivalent to the weak-type $\BMO$-regularity of the couple~$\left(X^{\frac r
2}, Y^{\frac r 2}\right)$.  By Proposition~\ref {xyxm} we have
$$
\left(X^{\frac r 2}\right)' Y^{\frac r 2} = \left(X^{\frac r 2}
Z^{\frac r 2}\right)' Z^{\frac r 2} Y^{\frac r 2} = \left(X^{\frac r 2} Z^{\frac
r 2}\right)' Y^{\frac r 2} Z^{\frac r 2},
$$
so condition~\ref {p:c:xywtbmor} is
equivalent to the weak-type $\BMO$-regularity of~$\left([X Z]^{\frac r 2}, [Y
Z]^{\frac r 2}\right)$, which is equivalent to condition~\ref {p:c:xzyzwtbmor} again
by Corollary~\ref {wtbmounif}.

With the help of the equivalence of conditions~\ref {p:c:xywtbmor} and~\ref
{p:c:xzyzwtbmor}, the equivalence of~\ref {p:c:xywtbmor} and~\ref
{p:c:yrxrwtbmor} follows in the standard way (see, e.~g., the proof of \cite [Theorem~8] {rutsky2011en}).
Condition~\ref {p:c:xywtbmor} is equivalent to the weak-type
$\BMO$-regularity of~$(X^r, Y^r)$, which is equivalent to the same of
$$
(X^r
(X^r)', Y^r (X^r)') = (\lclassg {1}, (X^r)' Y^r) = ((Y^r)' Y^r, (X^r)' Y^r),
$$
which is equivalent to condition~\ref {p:c:xzyzwtbmor}.

The symmetry, which is trivial for the usual $\BMO$-regularity, seems to require
the self-duality of the $\BMO$-regularity property (see \cite [Theorem~2]
{kisliakov2002en}, \cite [Theorem~1] {rutsky2011en} and~\cite [Theorem~1]
{rutsky2015en}).
By Proposition~\ref {wtbmoruniv} condition~\ref {p:c:xywtbmor} is equivalent to
the $\BMO$-regularity of lattice~$Z_1 = \left(\lclassg {1}, {\lclassg
{1}\strut}^{\frac 1 2} {\left[(X^r)' Y^r\right]\strut}^{\frac 1 2}\right)_{\theta, 1}$ with some~$0 < \theta < 1$,
which is equivalent to the $\BMO$-regularity of~$Z_1^{\frac 1 2} =
\left(\lclassg {2}, {\lclassg {2}\strut}^{\frac 1 2} {\left[(X^r)'^{\frac 1 2}
(Y^r)^{\frac 1 2}\right]\strut}^{\frac 1 2}\right)_{\theta, 2}$ by
Proposition~\ref {p:riscale}.
Since both lattices in the latter couple have order continuous norm, their
intersection is dense in each of them, and their Banach duals coincide with the
order duals.  Moreover, the intersection of these spaces is separable by~\cite
[Chapter~IV, \S 3, Theorem~3] {kantorovichold}, and it is also dense in the
interpolation space~$Z_1^{\frac 1 2}$ by~\cite [Theorem~3.4.2 (b)] {bergh}, so
by the same theorem from~\cite {kantorovichold} lattice~$Z_1^{\frac 1 2}$ has order continuous norm. 
Therefore, by the duality theorem for the real interpolation~\cite
[Theorem~3.7.1] {bergh} and~\cite [Theorem~1] {rutsky2011en} the
$\BMO$-regularity of~$Z_1$ is equivalent to the $\BMO$-regularity of the dual
lattice $$ \left(Z_1^{\frac 1 2}\right)' = \left(\lclassg {2}', \left[{\lclassg {2}\strut}^{\frac 1 2} {\left[(X^r)'^{\frac 1 2} (Y^r)^{\frac 1 2}\right]\strut}^{\frac 1 2}\right]'\right)_{\theta, 2} = \left(\lclassg {2},
{\lclassg {2}\strut}^{\frac 1 2} {\left[(Y^r)'^{\frac 1 2}
(X^r)^{\frac 1 2}\right]\strut}^{\frac 1 2}\right)_{\theta, 2}.
$$
Raising it to the power~$2$ with the help of Proposition~\ref {p:riscale} and making
use of Proposition~\ref {wtbmoruniv} yields the equivalence to the weak-type
$\BMO$-regularity of~$(Y, X)$,
which is condition~\ref {p:c:yxwtbmor}.

\section {Bounded AK-stability and proof of Theorem~\ref {akbmoreq0}}
\label {basS}

We need the following stronger species of the $\AK$-stability
property.  They were introduced in~\cite {kisliakov2003} and in~\cite
{rutsky2010en} respectively, but appeared implicitly in earlier research.
\begin {definition}
\label {sakstability}
A quasi-normed couple~$(X, Y)$ of lattices of measurable functions on~$\mathbb T
\times \Omega$ is called strongly $\AK$-stable with constant~$C$
if for any~$H \in (X + Y)_A$
and $f \in X$, $g \in Y$ such that $H = f + g$ there exist some $F \in X_A$, $G
\in Y_A$ such that~$H = F + G$ and~$\|F\|_X \leqslant C \|f\|_X$,
$\|G\|_Y \leqslant C \|g\|_Y$.
\end {definition}

The distinction between the $\AK$-stability and the strong $\AK$-stability
appears to be mostly technical in nature; see Proposition~\ref
{saksteq} below and remarks before it.

\begin {definition}
\label {bakstability}
A quasi-normed couple~$(X, Y)$ of lattices of measurable functions on~$\mathbb T
\times \Omega$ is called boundedly $\AK$-stable with constant~$C$ if for any $f
\in X$ and $g \in Y$ there exists some $U \in \hclass {\infty} {\mathbb T
\times \Omega}$ such that~$\|U\|_{\hclassg {\infty}} \leqslant C$, $\|U g\|_X
\leqslant C \|f\|_X$ and $\|(1 - U) f\|_Y \leqslant C \|g\|_Y$.
\end {definition}

The meaning of the bounded $\AK$-stability is clear from the following
reformulation (clarifying the discussion in~\cite [Section~1.3] {rutsky2010en}),
which also easily generalizes to the $\AK$-stability of several lattices.
\begin {proposition}
\label {bakstx}
Suppose that~$(X, Y)$ is a couple of quasi-normed lattices of measurable
functions on~$\mathbb T \times \Omega$ satisfying property~$(*)$.  Then it is
boundedly $\AK$-stable if and only if for any~$H \in (X + Y)_A$ and~$f \in
X$, $g \in Y$ such that~$H = f + g$ there exist some~$U \in \hclassg
{\infty}$ such that~$\|U H\|_X \leqslant c \|f\|_X$, $\|(1 - U) H\|_Y \leqslant c \|g\|_Y$
and~$\|U\|_{\hclassg {\infty}} \leqslant c$ with a constant~$c$ independent
of~$f$, $g$ and~$H$.
\end{proposition}
Indeed, it is easy to see that the bounded $\AK$-stability of a couple~$(X, Y)$
implies its strong $\AK$-stability with decompositions of the form~$F = U H$ and $G = (1
- U) H$ (see, e.~g., the proof of~\cite [Proposition~4] {rutsky2010en}).
Conversely, suppose that we are given some~$f \in X$ and~$g \in Y$.  We may
assume that these functions are nonnegative and, moreover, $\log f (\cdot,
\omega), \log g (\cdot, \omega) \in \lclassg {1}$ for almost all~$\omega \in \Omega$ by
propety~$(*)$.
We construct an outer function~$H = \exp \left(\log [f + g] + i \mathcal H [f + g]\right)$ such that~$|H| = f + g$ almost everywhere.
Here~$\mathcal H$ denotes the Hilbert transform acting in the first variable. 
Then function~$U$ from the assumptions also satisfies Definition~\ref
{bakstability} since~$|U| g \leqslant |U H|$ and~$|1 - U| f \leqslant |(1 - U)
H|$ almost everywhere.
  
We will see that for couples satisfying the
assumptions of either Theorem~\ref {akbmoreq0} or Theorem~\ref {goodrealinterp}
the $\AK$-stability is equivalent to the bounded $\AK$-stability.
It is easy to establish this equivalence in the important case $X = \lclassg
{\infty}$.
\begin {proposition} {\cite [Proposition~3] {rutsky2010en}}
\label {bakstlinf}
Let~$Z$ be a quasi-normed lattice of measurable functions on~$\mathbb T \times
\Omega$ satisfying property~$(*)$.  Couple~$\left(\lclassg {\infty}, Z\right)$
is strongly $\AK$-stable if and only if it is boundedly $\AK$-stable.
\end{proposition}

\begin {proposition} {\cite [Proposition~2] {rutsky2010en}}
\label {bakmult}
Let~$X, Y$ and~$Z$ be quasi-normed lattices of measurable functions on~$\mathbb
T \times \Omega$ such that~$(X, Y)$ is boundedly $\AK$-stable.  Then $(X Z, Y
Z)$ is also boundely $\AK$-stable.
\end{proposition}

These two simple results imply that a $\BMO$-regular couple~$(X, Y)$ is
boundedly $\AK$-stable, which was already noted in~\cite [\S 1.3] {rutsky2010en}.
For clarity, let us spell out an argument proving this.
For a weight~$\weightw$, which for simplicity we assume to be positive almost
everywhere, and for a lattice~$Z$ the weighted lattice~$Z (\weightw)$ is defined
by $Z (\weightw) = \{\weightw f \mid f \in Z\}$ with norm~$\|g\|_{Z (\weightw)} =
\left\|g \weightw^{-1}\right\|_Z$.
If~$\weightu$
and~$\weightv$ are some $\BMO$-majorants in the sense of Definition~\ref {bmorcdef} then $(\lclass
{\infty} {\weightu}, \lclass {\infty} {\weightv})$ is boundedly $\AK$-stable,
which was already implicit in the proof at the end of~\cite [\S 3.4]
{kisliakov1999}.  But for this couple the bounded $\AK$-stability
follows from the usual one: $\left(\lclassg {\infty}, \lclass {\infty}
{\weightu^{-1} \weightv}\right)$ is $\AK$-stable, so by Proposition~\ref
{bakstlinf} it is boundedly $\AK$-stable, and we may apply Proposition~\ref {bakmult} with~$Z = \lclass
{\infty} {\weightu}$.

The following result is a substantial improvement over~\cite [Proposition~6]
{rutsky2010en}.
\begin {proposition}
\label {akbakxprimey}
Let~$X$ and~$Y$ be Banach lattices of measurable functions on~$\mathbb T \times
\Omega$ satisfying the Fatou property and property~$(*)$.
Suppose also that~$X' Y$ is a Banach space.
Then~$(X, Y)$ is $\AK$-stable if and only if it is
boundedly~$\AK$-stable.
\end {proposition}
Indeed, if~$(X, Y)$ is $\AK$-stable then by~\cite [Lemma~4] {kisliakov2003} so
is~$(X' X, X' Y) = (\lclassg {1}, X' Y)$, by~\cite [Lemma~7]
{kisliakov2002en} couple~$(\lclassg {1}', (X' Y)') = (\lclassg {\infty}, (X'
Y)')$ is $\AK$-stable, and by Proposition~\ref {bakstlinf} it is boundedly
$\AK$-stable.  By Proposition~\ref {bakmult} and
Proposition~\ref {xyxm} couple~$(\lclassg {\infty} Y, (X' Y)' Y) = (Y, (X')') =
(Y, X)$ is then boundedly $\AK$-stable.

The following observation is well known.  It is an easy consequence of the
definition of the Calder\'on-Lozanovski\u{\i} product and the Young
inequality (see, e.~g., the proof of \cite [Lemma~5] {kisliakov2003}).
\begin {proposition}
\label {czprodctype}
Let~$X_0$ and~$X_1$ be some quasi-normed lattices of measurable functions
on the same measurable space.  Then the Calder\'on-Lozanovski\u{\i} product
$X_0^{1 - \theta} X_1^\theta$ is a space of type~$\mathcal C_\theta (X_0, X_1)$ for all~$0 < \theta < 1$.  That is, $(X_0, X_1)_{\theta, 1} \subset X_0^{1
- \theta} X_1^\theta \subset (X_0, X_1)_{\theta, \infty}$.
\end{proposition}

The following observation is key for the ``if'' part of Theorem~\ref
{akbmoreq0}.
\begin {proposition}
\label {linfbmorbakst}
Suppose that~$Z$ is a quasi-normed lattice of measurable functions on~$\mathbb T
\times \Omega$ such that~$(\lclassg {\infty}, Z)_{\theta, \infty}$ is
$\BMO$-regular with some~$0 < \theta < 1$.  Then~$(\lclassg {\infty}, Z^\alpha)$
is boundedly $\AK$-stable for all~$0 < \alpha < \theta$.
\end{proposition}
Indeed, let~$f \in \lclassg {\infty}$ and~$g \in Z^\alpha$.
For simplicity we may
assume that~$g$ is nonnegative and~$f = 1$. 
Then~$g^{\frac \theta \alpha} \in Z^\theta \subset (\lclassg {\infty},
Z)_{\theta, \infty}$ with norm at most~$c C \|g\|_{Z^\alpha}^{\frac \theta
\alpha}$  by
Proposition~\ref {czprodctype}, where~$c$ is some constant independent of~$g$.
Lattice~$\lclassg {\infty}$ is $\BMO$-regular, so couple~$\left(\lclassg
{\infty}, (\lclassg {\infty}, Z)_{\theta, \infty}\right)$ is also
$\BMO$-regular, and hence it is boundedly $\AK$-stable.
Thus there exists some~$U \in \hclassg {\infty}$ such that~$\|U\|_{\lclassg
{\infty}} \leqslant C$, $\left\| U
g^{\frac \theta \alpha}\right\|_{\lclassg {\infty}} \leqslant C \|f\|_{\lclassg {\infty}} = C$ and~$\| (1 - U) f\|_{(\lclassg {\infty},
Z)_{\theta, \infty}} \leqslant c C \|g\|_{Z^\alpha}^{\frac \theta \alpha}$.
From the second estimate it follows that
\begin {equation}
\label {eqbkst1}
\left\|U g\right\|_{\lclassg {\infty}} =
\left\||U|^{\frac \theta \alpha} g^{\frac \theta \alpha}\right\|_{\lclassg
{\infty}}^{\frac \alpha \theta} \leqslant C^{1 - \frac \alpha \theta} \left\|U
g^{\frac \theta \alpha}\right\|_{\lclassg {\infty}}^{\frac \alpha \theta}
\leqslant C.
\end {equation}
But we also have~$\| (1 - U) f \|_{\lclassg {\infty}} \leqslant C + 1$.
By Proposition~\ref {czprodctype} and the reiteration theorem  lattice~$Z^\alpha$ is a space of
type~$\mathcal C_{\frac \alpha \theta} \left(\lclassg
{\infty}, (\lclassg {\infty}, Z)_{\theta, \infty}\right)$, thus
\begin {equation}
\label {eqbkst2}
\| (1 - U) f \|_{Z^{\alpha}} \leqslant \| (1 - U) f \|_{\lclassg
{\infty}}^{1 - \frac \alpha \theta} \| (1 - U) f \|_{(\lclassg {\infty},
Z)_{\theta, \infty}}^{\frac \alpha \theta} \leqslant
c_1 \|g\|_{Z^\alpha}
\end {equation}
with some constant~$c_1$ independent of~$g$.  \eqref {eqbkst1} and~\eqref
{eqbkst2} together show that couple~$(\lclassg {\infty}, Z^\alpha)$ is indeed
boundely $\AK$-stable.

\begin {lemma}
\label {bakstabilitygen}
Suppose that~$\alpha, \beta > 0$. A couple~$(X, Y)$ of quasi-normed lattices of
measurable functions on~$\mathbb T \times \Omega$ is boundedly $\AK$-stable with
some constant~$C$ if and only if for any~$f \in X$ and~$g \in Y$ there exists
some $V \in \hclass {\infty} {\mathbb T \times \Omega}$ such that~$\|V\|_{\hclassg
{\infty}} \leqslant C'$, $\left\||V|^\alpha g\right\|_X \leqslant C' \|f\|_X$
and $\left\||1 - V|^\beta f\right\|_Y \leqslant C' \|g\|_Y$ with some
constant~$C'$.
\end{lemma}
To prove the ``if'' part, take some integer numbers~$M \geqslant \alpha$ and~$N \geqslant \beta$,
let~$f \in X$ and~$g \in Y$ be some nonnegative functions, and a
corresponding function~$V$ from the satement of the lemma.
Functions~$V^M$ and~$(1 - V)^N$ satisfy the
assumptions of the corona theorem (see, e.~g., \cite [Proposition~2] {kisliakovrutsky2012en}), so there exist
some~$U_0, U_1 \in \hclassg {\infty}$ such that~$V^M U_0 + (1 - V)^N U_1 = 1$
and~$\|U_j\|_{\hclassg {\infty}} \leqslant c_1$, $j \in \{0, 1\}$ with
some~$c_1$ independent of~$f$ and~$g$.
Let~$U = V^M U_0$.  Then~$|U| \leqslant c_1 |V|^{M} \leqslant c_1
C'^{M - \alpha} |V|^\alpha$, and similarly~$|1 - U| \leqslant c_1 C'^{N - \beta} |1 - U|^\beta$,
which yields the claimed estimates.
For the ``only if'' part we choose some integer numbers~$M \geqslant \frac 1
\alpha$, $N \geqslant \frac 1 \beta$, apply the corona theorem to find some
bounded analytic functions~$V_0$ and~$V_1$ satisfying~$U^M V_0 + (1 - U)^N V_1 = 1$ with
suitable estimates, and take $V = U^M V_0$. 

We mention that the use of the corona theorem can be easily avoided in the proof
of Lemma~\ref {bakstabilitygen}, perhaps at a slight expense of clarity and
generalizations to several lattices.
In the ``if'' part we may first take~$V_1 = V^M$, which yields estimates~$\left\||V| g\right\|_X \leqslant C'' \|f\|_X$
and~$\left\||1 - V|^\beta f\right\|_Y \leqslant C'' \|g\|_Y$ from the
assumptions with some suitable constant~$C''$, since~$|1 - V_1| = \left|(1 - V)
\sum_{j = 0}^{M - 1} V^j \right| \leqslant M C'^M |1 - V|$.  Taking~$U = 1 - (1
- V_1)^N$ then yields the bounded $\AK$-stability of~$(X, Y)$ by a similar
estimate, and the ``only if'' part is treated in the same way.

Observe that by the homogeneity the conditions of Lemma~\ref
{bakstabilitygen} may be restricted to~$\|g\|_Y = 1$.  Further replacing~$f$
with~$f_1$ such that~$f = s f_1$ and~$s = \|f\|_X$ yields the following characterization of bounded $\AK$-stability
that will be used in the proof of Theorem~\ref {t:erbakst} below.
\begin {corollary}
\label {c:tbakstabilitygen}
Suppose that~$\alpha, \beta > 0$. A couple~$(X, Y)$ of quasi-normed lattices of
measurable functions on~$\mathbb T \times \Omega$ is boundedly $\AK$-stable with
some constant~$C$ if and only if for any~$s > 0$ and~$f \in X$, $g \in Y$ such
that~$\|f\|_X = \|g\|_Y = 1$ there exists some~$V \in \hclass {\infty} {\mathbb
T \times \Omega}$ such that~$\|V\|_{\hclassg {\infty}} \leqslant C'$,
$\left\||V|^\alpha g\right\|_X \leqslant C' s$ and~$\left\||1 - V|^\beta
f\right\|_Y \leqslant C' s^{-1}$ with some constant~$C'$.
\end{corollary}

The following observation is a generalization of~\cite [Proposition~1]
{rutsky2010en}, where the case~$\delta < 1$ was trivially established.
See also~\cite [Theorem~2] {kisliakov2003} and~\cite [Theorem~3.6]
{kisliakov1999}.
It is interesting to note that the latter theorem together with~\cite
[Proposition~1] {rutsky2010en} already implies that if a couple~$(X, Y)$ of
quasi-normed lattices is boundedly $\AK$-stable then the couple~$(X^\delta,
Y^\delta)$ is $\AK$-stable for all~$\delta > 0$.
However, it does not seem to say anything about the bounded $\AK$-stability, 
and both~\cite [Theorem~2] {kisliakov2003} and~\cite [Theorem~3.6]
{kisliakov1999} are rather nontrivial results. 
\begin {proposition}
\label {bakrscale}
Suppose that~$(X, Y)$ is a couple of quasi-normed lattices of measurable
functions on~$\mathbb T \times \Omega$.  If~$(X, Y)$ is boundedly $\AK$-stable
then so is~$(X^\delta, Y^\delta)$ for all~$\delta > 0$.
\end{proposition}
Indeed, if~$f \in
X^\delta$ and~$g \in X^\delta$ then by the bounded $\AK$-stability of~$(X, Y)$
there exists some~$V \in \hclassg {\infty}$ such that
$$
\|V\|_{\lclassg {\infty}} \leqslant C, \quad \left\|V g^{\frac 1
\delta}\right\|_X \leqslant C \left\|f^{\frac 1 \delta}\right\|_X \text { and }
\left\|(1 - V) f^{\frac 1 \delta} \right\|_Y \leqslant C \left\|g^{\frac
1 \delta}\right\|_Y.
$$
These conditions are exactly~$\left\| |V|^\delta 
g\right\|_{X^\delta} \leqslant C^\delta \|f\|_{X^\delta}$ and~$\left\| |1 - V|^\delta f \right\|_{Y^\delta} \leqslant
C^\delta \|g\|_{Y^\delta}$, so by Lemma~\ref {bakstabilitygen}
couple~$(X^\delta, Y^\delta)$ is boundedly $\AK$-stable.

We are now ready to prove the ``if'' part of Theorem~\ref {akbmoreq0} under the
assumption that lattice~$X' Y$ is Banach.
\begin {proposition}
\label {akbmoreq0i}
Let~$(X, Y)$ be a couple of quasi-Banach lattices of measurable
functions on~$\mathbb T \times \Omega$ that are $r$-convex with some $r > 0$
satisfying the Fatou property and property~$(*)$. Suppose also that~$(X^r)'
Y^r$ is a Banach lattice.  If~$(X, Y)$ is weak-type $\BMO$-regular then it is boundedly $\AK$-stable.
\end{proposition}
By Proposition~\ref {wtbmoruniv} lattice
$
\left(\lclassg {1},
\lclassg {1}^{1 - \beta} \left[(X^r)' (Y^r)\right]^\beta\right)_{\eta, 1}
$
is $\BMO$-regular
with some $0 < \beta, \eta < 1$ such
that~$\beta \eta = \theta$.  Similarly to the proof of symmetry in
Proposition~\ref {wtbmordiv}, by the duality theorems for the real interpolation
and for $\BMO$-regularity \cite [Theorem~1] {rutsky2011en} lattice
$$
\left(\lclassg {1}, \lclassg {1}^{1 - \beta} \left[(X^r)'
(Y^r)\right]^\beta\right)_{\eta, 1}' = \left(\lclassg {\infty}, \left[(X^r)'
(Y^r)\right]'^\beta\right)_{\eta, \infty}
$$
is $\BMO$-regular.  By Proposition~\ref {linfbmorbakst} it follows that
couple~$\left(\lclassg {\infty}, \left[(X^r)'(Y^r)\right]'^\alpha\right)$ 
is boundedly $\AK$-stable for all~$0 < \alpha < \theta$.
Therefore, couple
$$
\left(\lclassg {\infty} Y^{r \alpha}, \left(\left[(X^r)'(Y^r)\right]'
Y^r\right)^\alpha\right) = \left(Y^{r \alpha},
\left[(X^r)'\right]'^\alpha\right) =
\left(Y^{r \alpha}, X^{r \alpha}\right)
$$
is also boundedly $\AK$-stable by Proposition~\ref {bakmult} with~$Z = Y^{r
\alpha}$.  Here we used the formula from Proposition~\ref {xyxm}.
Finally, by Proposition~\ref {bakrscale}
couple~$(X, Y)$ is boundedly $\AK$-stable.

The proof of the ``only if'' part of Theorem~\ref {akbmoreq0} is based on the
following result.\footnote {We already mentioned in the introduction that the
main results of the author published in~\cite {rutsky2013ben}, unfortunately, are
flawed.  The cited corollary and other
cited results from~\cite {rutsky2013ben}, however, are sound.}
\begin {proposition} {\cite [Corollary~13] {rutsky2013ben}}
\label {beq0corr}
Let~$Z$ be a Banach lattice of measurable functions on~$\mathbb T \times \Omega$
satisfying the Fatou property and property~$(*)$, and assume that lattices $Z$
and~$Z'$ have order continuous norm.
Suppose that the couple~$(Z, Z')$ is $\AK$-stable.
Then the Riesz projection is bounded in~$(\lclassg {2}, Z)_{\zeta, 2}$
for all sufficiently small~$0 < \zeta < 1$.
\end {proposition}
We establish the ``only if'' part of Theorem~\ref {akbmoreq0} under a somewhat
broader set of assumptions,
which cover the same assumptions as in the main results of~\cite
{rutsky2013ben} and follow the corresponding details of the reductions.
However, we generalize them substantially, and at least some further
generalizations appear to be possible.
Specifically, \cite {rutsky2013ben} only had assumption~\ref {c:p:convprime} below with~$p = 2$,
assumption~\ref {c:p:lp} with~$p \in \{1, 2, \infty\}$, and it did not have
assumptions~\ref {c:p:convconc} and~\ref {c:p:convconv}.  From assumption~\ref
{c:p:lp} it follows that if one of the lattices is simultaneously $p$-convex and $p$-concave then we
need no restrictions on the other lattice.
For simplicity, we state these
assumptions with some asymmetry and redundancy.  Most prominently,
assumption~\ref {c:p:convconv} generalizes assumptions~\ref {c:p:convprime},
\ref {c:p:convconc} and~\ref {c:p:lp}.
\begin {proposition}
\label {akbmoreq0r}
Suppose that~$(X, Y)$ be a couple of Banach lattices of
measurable functions on~$\mathbb T \times \Omega$ satisfying
the Fatou property, property~$(*)$ and at least one of the following
conditions:
\begin {enumerate} [label=\text {\upshape(\roman*)}, leftmargin=*, widest=iii]
  \item \label {c:p:xxprime} Lattices~$X$ and~$Y$ have order continuous norm
  and~$Y = X'$;
  \item \label {c:p:convprime} $X$ is $p$-convex and~$Y$ is $p'$-convex with
  some~$1 \leqslant p \leqslant \infty$;
  \item \label {c:p:xprimey} $X' Y$ is Banach;
  \item \label {c:p:convconc} $X$ is $p$-concave and~$Y$ is $p$-convex with
  some~$1 \leqslant p \leqslant \infty$;
  \item \label {c:p:lp} $X = \lclassg {p}$ with~$1 \leqslant p \leqslant
  \infty$;
  \item \label {c:p:convconv} $X$ is $p$-convex and $q$-concave with some $1
  \leqslant p \leqslant q \leqslant \infty$ and~$Y$ is $\left(\frac 1 {p'}
  + \frac 1 q\right)^{-1}$-convex.
\end{enumerate}
If couple~$(X, Y)$ is $\AK$-stable then it is weak-type $\BMO$-regular.
\end{proposition}

Suppose that under the
assumptions of Proposition~\ref {akbmoreq0r} couple~$(X, Y)$ is $\AK$-stable.
Under assumptions~\ref {c:p:xxprime}, Proposition~\ref {beq0corr} applied to~$Z
= Y$ directly shows that the Riesz projection is
bounded in~$Z_1 = (\lclassg {2}, Y)_{\zeta, 2}$ with some~$0 < \zeta < 1$, which
implies by~\cite [Theorem~3] {kisliakov2002en} (see also~\cite {rutsky2011en})
that~$Z_1$ is $\BMO$-regular.  Lattice~$\left(\lclassg {2},
Y\right)_{\theta, 2}$ is $\BMO$-regular with some~$0 < \theta < 1$.This yields
the $\BMO$-regularity of~$\left(\lclassg {2}, Y\right)_{\theta, 2}^2 =
\left(\lclassg {1}, Y^2\right)_{\theta, 1} = \left(\lclassg {1},
X' Y\right)_{\theta, 1}$, which is the weak-type $\BMO$-regularity of~$(X, Y)$.

First we will show that the conclusion holds true if the assumptions~\ref
{c:p:convprime} are satisfied in the special case~$p = 2$.
Let~$X_1 = X^2$ and~$Y_1 = Y^2$.
Thus~$(X, Y) = \left(X_1^{\frac 1 2}, Y_1^{\frac 1
2}\right)$ is an $\AK$-stable couple of Banach lattices.
Let~$Z = {\lclassg {2}\strut}^{\frac
1 2} {X_1'\strut}^{\frac 1 4} {Y_1^{\phantom '}\strut}^{\frac 1 4}$.  A
simple computation in the proof of \cite [Theorem~2] {rutsky2013ben} shows
that couple~$(Z, Z')$ is $\AK$-stable, and it satisfies assumptions~\ref
{c:p:xxprime}.
Therefore, lattice
$$
Z_1^2 = \left(\lclassg {1}, \lclassg {1}^{\frac 1 2}
\left({X_1\strut}' {Y_1\strut}\right)^{\frac 1 2}\right)_{\zeta, 1} =
\left(\lclassg {1}, {X_1\strut}' {Y_1\strut}\right)_{\frac \zeta 2, 1} =
\left(\lclassg {1}, (X^2)' Y^2\right)_{\frac \zeta 2, 1}
$$
is also $\BMO$-regular, so couple~$(X, Y)$ is weak-type $\BMO$-regular by
Proposition~\ref {wtbmoruniv} as claimed.  Here we have used Proposition~\ref
{p:riscale} and the reiteration formula.

If assumptions~\ref {c:p:xprimey} are satisfied, by Proposition~\ref
{akbakxprimey} couple~$(X, Y)$ is boundedly $\AK$-stable, and by Proposition~\ref {bakrscale}
so is couple~$(X^{\frac 1 2}, Y^{\frac 1 2})$.  This couple satisfies
assumptions~\ref {c:p:convprime} with~$p = 2$, so it is weak-type
$\BMO$-regular, and Proposition~\ref {wtbmoruniv} again yields the weak-type $\BMO$-regularity of the original couple~$(X, Y)$.

Under assumptions~\ref {c:p:convconc} lattice~$X'$ is $p'$-convex, so
lattice~$X' Y$ is $1$-convex and the couple satisfies assumptions~\ref
{c:p:xprimey}.

Under assumptions~\ref {c:p:lp} cases~$p = 1$ and~$p = \infty$ satisfy
assumptions~\ref {c:p:xprimey} (in the case~$p = \infty$ we need to reverse the
order of the couple), so the interesting case is~$1 < p < \infty$.  We may further assume
that $p \geqslant 2$, otherwise we may pass to the duals in the $\AK$-stability
by \cite [Lemma~7] {kisliakov2002en} and then use the duality in
Proposition~\ref {wtbmordiv}.
The $\AK$-stability of $(\lclassg {p}, Y) =
\left( Y^{\frac 1 p} \left(Y'^{\frac {p'} p}\right)^{\frac 1 {p'}},
{Y\strut}^{\frac 1 p} {Y\strut}^{\frac 1 {p'}}\right)$ by~\cite [Corollary to
Lemma~4] {kisliakov2003} is equivalent to the $\AK$-stability of~$\left( \left(Y'^{\frac {p'} p}\right)^{\frac 1 {p'}},
Y^{\frac 1 {p'}}\right) = \left( Y'^{\frac 1 p},
Y^{\frac 1 {p'}}\right)$.
Let~$Z_1 = Y'^{\frac 1 2 \left(\frac 1 {p'} -
\frac 1 p\right)} \lclassg {p}^{\frac 1 2}$.  By~\cite [Lemma~4]
{kisliakov2003} the latter $\AK$-stability implies the $\AK$-stability of
\begin {multline*}
\left( Y'^{\frac 1 p} Z_1, Y^{\frac 1 {p'}} Z_1\right) =
\\
\left( Y'^{\frac 1 p + \frac 1 2 \left(\frac 1 {p'} - \frac 1 p\right)}
\lclassg {p}^{\frac 1 2}, Y^{\frac 1 {p'} - \frac 1 2 \left(\frac 1 {p'}
- \frac 1 p \right)} Y^{\frac 1 2 \left(\frac 1 {p'}
- \frac 1 p \right)} Y'^{\frac 1 2 \left(\frac 1 {p'}
- \frac 1 p \right)} \lclassg {p}^{\frac 1 2}\right)
=
\\
\left({Y'\strut}^{\frac 1 2} {\lclassg {p}\strut}^{\frac 1 2}, {Y\strut}^{\frac
1 2}
{\lclassg {1}\strut}^{\frac 1 2 \left(\frac 1 {p'}
- \frac 1 p \right)} {\lclassg {p}\strut}^{\frac 1 {2 p}}\right)
=
\left({Y'\strut}^{\frac 1 2} {\lclassg {p}\strut}^{\frac 1 2}, {Y\strut}^{\frac
1 2} {\lclassg {p'}\strut}^{\frac 1 2}\right)
=
(Z, Z')
\end {multline*}
with~$Z = {Y'\strut}^{\frac 1 2} {\lclassg {p}\strut}^\frac 1 2$.
This couple satisfies assumptions~\ref {c:p:xxprime}, so it is weak-type
$\BMO$-regular.
By running the respective multiplications and divisions in reverse we get
the weak-type $\BMO$-regularity of the original couple~$(\lclassg {p}, Y)$
by Proposition~\ref {wtbmordiv}. 

Now, suppose that assumptions~\ref {c:p:convconv} are satisfied.  If~$p =
\infty$ then assumptions~\ref {c:p:convconc} are satisfied for couple~$(Y, X)$,
so the interesting case is~$1 \leqslant p < \infty$. Let~$Z_2 = (X^p)'^{\frac 1
p}$.
Observe that~$\lclassg {1} = X^p \left(X^p\right)'$, so~$\lclassg {p} = \lclassg
{1}^{\frac 1 p} = X Z_2$. Couple $\left(X Z_2, Y Z_2\right) = \left(\lclassg
{p}, Y Z_2\right)$ is $\AK$-stable by~\cite [Lemma~4] {kisliakov2003}.  Lattice~$X^p$ is $\frac q p$-concave,
so~$(X^p)'$ is $\left(\frac q p\right)'$-convex and~$Z$
is $r$-convex with~$r = p \left(\frac q p\right)' = \frac {q p} {q - p}$.
Simple computations show that~$Z Y$ is then $1$-convex.  Thus, this couple
satisfies assumptions~\ref {c:p:lp}, so it is weak-type $\BMO$-regular,
and~$(X, Y)$ is weak-type $\BMO$-regular by Proposition~\ref {wtbmordiv}.
Finally, assumptions~\ref {c:p:convprime} with arbitrary~$p$ imply
assumptions~\ref {c:p:convconv} with~$q = \infty$.

\begin {corollary}
\label {akbmoreq0rb}
Suppose that~$(X, Y)$ is a couple of $r$-convex quasi-normed lattices of
measurable functions on~$\mathbb T \times \Omega$ with some~$r > 0$ satisfying
the Fatou property and property~$(*)$.  If~$(X, Y)$ is boundedly $\AK$-stable
then~$(X, Y)$ is weak-type $\BMO$-regular.
\end {corollary}
Indeed, by Proposition~\ref {bakrscale} couple~$\left(X^{\frac r 2},
Y^{\frac r 2}\right)$ is boundedly $\AK$-stable, so by Proposition~\ref
{akbmoreq0r} it is weak-type $\BMO$-regular, and by Corollary~\ref
{wtbmounif} couple~$(X, Y)$ is weak-type $\BMO$-regular.

\section {Bounded $\nplus$-stability}

\label {akbakeqs}

Similarly to the bounded $\AK$-stability, we may define a stronger version of
$\nplus$-stability with respect to the real interpolation functors. 
Let~$\lambda > 1$.  We may fix the standard value~$\lambda = 2$ in the present
work.
Recall that~$J (t, g; X, Y) = \|g\|_X \vee t \|g\|_Y$ for~$t
> 0$ and~$g \in X \cap Y$, and the real interpolation space~$(X, Y)_{\theta, p}$
may be defined by the so-called~$J$ method as the space of functions~$f \in X +
Y$ having decompositions~$f = \sum_j f_j$ with finite norm
$\left\|\{f_j\}_{j \in \mathbb Z}\right\|_{\mathcal X_{(X, Y)_{\theta, p}}} =
\left\|\left\{\lambda^{-\theta j} J (\lambda^j, f_j; X, Y) \right\}_{j \in \mathbb Z}\right\|_{\lsclass {p}}$, and the norm
of~$f$ in~$(X, Y)_{\theta, p}$ is taken to be the infinum of~$\left\|\{f_j\}_{j
\in \mathbb Z}\right\|_{\mathcal X_{(X, Y)_{\theta, p}}}$ over all such
decompositions.
\begin {definition}
\label {bnplusstdef}
Let~$(X, Y)$ be a couple of quasi-Banach lattices of measurable functions
on~$\mathbb T \times \Omega$, and let~$\lambda > 1$.  We say
that~$(X, Y)$ is boundedly $\nplus$-stable with respect to~$(\cdot, \cdot)_{\theta, p}$ with constant~$C$
if for any~$f \in \left[(X, Y)_{\theta, p}\right]_A$ there exists
some~$\varphi = \left\{\varphi_j\right\}_{j \in \mathbb Z} \in \hclass {\infty}
{\lsclass {1}}$ with norm at most~$C$ such that
$\sum_j \varphi_j = 1$ and $\left\|\left\{\lambda^{-\theta j} J (\lambda^j, 
\varphi_j f; X, Y) \right\}_{j \in \mathbb Z}\right\|_{\lsclass {p}} \leqslant C
\|f\|_{(X, Y)_{\theta, p}}$.
\end{definition}
It is easy to verify that this definition does not depend on a particular choice
of the parameter~$\lambda > 1$.
Unlike the bounded $\AK$-stability, it is not
clear if this property is stable with respect to the multiplication by a lattice as in Proposition~\ref
{bakmult}.
Also, it is not clear if there are easy and general equivalence results deriving it from the usual $\nplus$-stability
similarly to Propositions~\ref {bakstlinf} and~\ref {akbakxprimey}, or if it is
stable under raising to powers greater than~$1$.  However, raising to powers~$0
< \delta < 1$ still works, and it allows us to improve convexity of boundedly
$\nplus$-stable lattices.
\begin {proposition}
\label {p:bnplusunif}
Suppose that a couple~$(X, Y)$ of quasi-Banach lattices of measurable
functions on~$\mathbb T \times \Omega$ satisfying property~$(*)$ is boundedly
$\nplus$-stable with respect to~$(\cdot, \cdot)_{\theta, p}$. 
Then~$\left(X^\delta, Y^\delta\right)$ is boundedly $\nplus$-stable with respect
to~$(\cdot, \cdot)_{\theta, \frac p \delta}$ for all~$0 < \delta < 1$.
\end{proposition}
Indeed, suppose that
$F \in \left[\left(X^\delta, Y^\delta\right)_{\theta, \frac p \delta}\right]_A$
with norm~$1$  under the assumptions of Propostion~\ref {p:bnplusunif}.
Lattice~$(X, Y)_{\theta, p}$ has property~$(*)$ by~\cite [Proposition~9]
{rutsky2013ben}, so there exists some~$g \geqslant |F|$, $\|g\|_{\left(X^\delta,
Y^\delta\right)_{\theta, \frac p \delta}} \leqslant 2$ such that~$\log g
(\cdot, \omega) \in \lclassg {1}$ for almost all~$\omega \in \Omega$.
We construct the corresponding outer function~$G = \exp \left(g + i
\mathcal H g\right)$.  Observe that~$G^{\frac 1 \delta} \in \left[(X,
Y)_{\theta, p}\right]_A$ with norm at most~$2^{\frac 1 \delta}$, so there exists
some~$\varphi = \left\{\varphi_j\right\}_{j \in \mathbb Z} \in \hclass {\infty}
{\lsclass {1}}$ with norm at most~$C$ satisfying Definition~\ref {bnplusstdef}
with~$f = G^{\frac 1 \delta}$.  Then the same function~$\varphi$ yields the
bounded $\nplus$-stability for~$F$ with respect to~$(\cdot, \cdot)_{\theta,
\frac p \delta}$, which follows from the estimate
$$
|\varphi_j F| \leqslant |\varphi_j G| = |\varphi_j|^{1 - \delta}
\left|\varphi_j G^{\frac 1 \delta}\right|^\delta \leqslant C^{1 - \delta}
\left|\varphi_j G^{\frac 1 \delta}\right|^\delta.
$$

With the
help of a fixed point theorem, we will show that bounded stability in the sense
of Definition~\ref {bnplusstdef} often naturally arises from the usual
stability.
We will do this in a more abstract setting that will be useful elsewhere.
\begin {definition}
\label {d:jxsp}
Let~$I \subset \mathbb Z$,
$\mathcal M$ be a $\sigma$-finite
measurable space.  Suppose that~$\mathcal X$ is a quasi-normed lattice of
measurable functions on~$\mathcal M \times I$ and~$R$ is a quasi-normed lattice
of measurable functions on~$\mathcal M$.
We say that~$\mathcal X$ is summable if for
any~$\{f_j\}_{j \in I} \in \mathcal X$ the sum~$\sum_{j \in I} |f_j|$ is finite
almost everywhere.
Let~$\mathcal S_n \{f_j\} = \sum_{j \in I \cap [-n, n]} f_j$, $n \in \mathbb
N \cup \{\infty\}$.
We say that~$\mathcal X$ is uniformly $R$-summable if~$\|\mathcal S_\infty -
\mathcal S_n\|_{\mathcal X \to R} \to 0$.
\end{definition}

\begin {definition}
\label {d:jxspaces}
Let~$\mathcal M$ be a $\sigma$-finite
measurable space, and suppose that~$\mathcal X$ is a summable quasi-normed
lattice of measurable functions on~$\mathcal M \times I$. 
We define a lattice
$$
J (\mathcal X) = \left\{\sum_{j \in I} f_j \mid
\{f_j\}_{j \in I} \in \mathcal X\right\}
$$
with the corresponding quasi-norm
$$\|f\|_{J (\mathcal X)} = \inf
\left\{\left\|\{f_j\}_{j \in I}\right\|_{\mathcal X} \mid \sum_{j \in I} f_j =
f, \{f_j\}_{j \in I} \in \mathcal X\right\}.
$$
\end{definition}
For example, if~$X$ and~$Y$ are quasi-Banach lattices,
taking~$I = \mathbb Z$
and a lattice~$\mathcal X_{(X, Y)_{\theta, p}}$ defined by
the norm
$$
\|\{f_j\}_{j \in \mathbb Z}\|_{\mathcal X_{(X, Y)_{\theta,
p}}} = \left\|\left\{\lambda^{-\theta j} J (\lambda^j, f_j; X, Y) \right\}_{j \in \mathbb Z}\right\|_{\lsclass {p}}
$$
yields~$J \left(\mathcal X_{(X, Y)_{\theta, p}}\right) = (X, Y)_{\theta,
p}$.  It is easy to see that~$\mathcal X_{(X, Y)_{\theta, p}}$ is uniformly $(X
+ Y)$-summable.
\begin {definition}
\label {d:jxstability}
Let~$I \subset \mathbb Z$ and let~$\mathcal X$ be a summable quasi-normed
lattice of measurable functions on~$\mathbb T \times \Omega \times I$.
We say that~$J (\mathcal X)$ is $\nplus$-stable with
constant~$C$ if for any~$f \in \left[J (\mathcal X)\right]_A$ there exists
some~$F = \{f_j\}_{j \in I} \in \mathcal X_A$ such that~$f = \sum_{j \in I} f_j$
and~$\|F\|_{\mathcal X} \leqslant C \|f\|_{J (\mathcal X)}$.
We say that~$J (\mathcal X)$ is boundedly $\nplus$-stable with constant~$C$ if
in the above we may take~$f_j = f \varphi_j$, $j \in I$, $\{\varphi_j\}_{j \in
\mathbb Z} \in \hclass {\infty} {\lsclass {1}}$ with norm at most~$C$.
\end{definition}
In the example above, the (bounded) $\nplus$-stability of~$\left(\mathcal X_{(X,
Y)_{\theta, p}}\right)$ is exactly the (bounded) $\nplus$-stability of the
couple~$(X, Y)$ with respect to~$(\cdot, \cdot)_{\theta, p}$.

\begin {theorem}
\label {t:jxstability}
Let~$I \subset \mathbb Z$ and let~$\mathcal X$ be a uniformly $R$-summable
Banach lattice of measurable functions on~$\mathbb T \times \Omega \times I$.
Suppose that~$\Omega$ is a discrete space, $\mathcal X$ has
the Fatou property and~$J (X)$ has property~$(*)$.
If~$J (\mathcal X)$ is $\nplus$-stable then it is boundedly
$\nplus$-stable.
\end{theorem}

The proof of Theorem~\ref {t:jxstability} is given in Section~\ref {s:jxsp} below. 
As a consequence of Proposition~\ref {p:bnplusunif}, we at once get the
following.
\begin {corollary}
\label {c:bnps}
Suppose that~$(X, Y)$ is a couple of Banach lattices of measurable functions
on~$\mathbb T \times \Omega$ with a discrete space~$\Omega$ satisfying the Fatou
property and property~$(*)$.  If~$(X, Y)$ is $\nplus$-stable with respect
to~$(\cdot, \cdot)_{\theta, p}$ then~$\left(X^\delta, Y^\delta\right)$ is
boundedly $\nplus$-stable with respect to~$(\cdot, \cdot)_{\theta, \frac p \delta}$ for
all~$0 < \delta \leqslant 1$.
\end{corollary}

The equivalence of $\AK$-stability and the bounded $\AK$-stability for Banach
lattices with the Fatou property and a discrete~$\Omega$ is an easy consequence
of Theorem~\ref {t:jxstability}.  The proof of Theorem~\ref {goodrealinterp}
also allows us to relax the convexity assumptions; see Corollary~\ref
{c:akbakeq} below.
We also mention that although we only work with
couples of lattices in this paper, the proof as written also yields the boundedness of the $\AK$-stability for an arbitrary finite family of
lattices.
\begin {corollary}
\label {akbakeq}
Suppose that~$X_0$, $X_1$ are Banach lattices of measurable functions
on~$\mathbb T \times \Omega$ satisfying the Fatou property and property~$(*)$.
Suppose also that~$\Omega$ is a discrete space.
Then~$(X_0, X_1)$ is $\AK$-stable if and only if it is boundedly~$\AK$-stable.
\end {corollary}
Let~$I = \{0, 1\}$.  Suppose that~$(X, Y)$ is $\AK$-stable. Let~$g_j \in X_j$
be some nonzero functions, and define a lattice~$\mathcal X$ on~$\mathbb T
\times \Omega \times I$ with the norm~${\|\{f_j\}_{j \in I}\|}_{\mathcal
X} = \bigvee_{j \in I} {\|g_j\|\strut}_{X_j}^{-1} {\|f_j\|\strut}_{X_j}$.
The $\AK$-stability implies the $\nplus$-stability of~$J (\mathcal X)$, and
incidentally is equivalent to the latter satisfied uniformly over
arbitrary functions~$g_j \in X_j$, $j \in I$.
Let~$h_j \geqslant |g_j|$, $j \in I$ be the corresponding majorants from
property~$(*)$, $h = \sum_{j \in I} |h_j|$ and~$H = \exp
\left(\log h + i \mathcal H \log h\right) \in \left[J (\mathcal X)\right]_A$
with norm at most~$2$.
By Theorem~\ref {t:jxstability} lattice~$J (\mathcal X)$ is boundedly
$\nplus$-stable with a constant~$C$, so there exist
some~$\left\{\varphi_j\right\}_{j \in I} \in \hclass {\infty} {\lsclass {1}}$
with norm at most~$C$ such that~$\sum_{j \in I} \varphi_j = 1$ and $$
\left\|\varphi_k \sum_{j \in I} |g_j|\right\|_{X_k} \leqslant \|\varphi_k
H\|_{X_k} \leqslant 2 C \|g_k\|_{X_k}, \quad k \in I.
$$
Setting~$U = \varphi_0$ yields the bounded $\AK$-stability of the couple~$(X_0,
X_1)$.

The uniform $R$-summability property is a convenient condition that ensures the
closedness with respect to the convergence in measure (on all sets of finite
measure) of the set of corresponding decompositions for lattices~$\mathcal X$ with the Fatou property,
and that lattices~$J (\mathcal X)$ inherit the Fatou property from~$\mathcal X$
\begin {proposition}
\label {p:jxsumcm}
Let~$(\mathcal M, \mu)$ be a $\sigma$-finite measurable space, $I \subset
\mathbb Z$, and let~$\mathcal X$ be a
uniformly $R$-summable Banach lattice of measurable functions on~$\mathcal M
\times I$ with the Fatou property.  Let~$C > 0$.
Then the graph of a set-valued mapping~$D_{\mathcal X, C} : J (\mathcal X) \to
2^{\mathcal X}$ defined by
$$
D_{\mathcal X, C} (f)
= \left\{\{f_j\}_{j \in I} \mid f = \sum_{j \in I} f_j, \|\{f_j\}_{j \in
I}\|_{\mathcal X} \leqslant C \right\}, \quad f \in J (\mathcal X)
$$
is closed with respect to the convergence in measure on sets of finite measure.
In particular, sets~$D_{\mathcal X, C} (f)$ are closed with respect to this
convergence, they are nonempty for~$C \geqslant \|f\|_{J (\mathcal X)}$ and $J
(\mathcal X)$ is a Banach lattice with the Fatou property.
\end{proposition}
Indeed, let~$f_k \in J (\mathcal X)$ and~$g^{(k)} \in D_{\mathcal X, C} (f_k)$
be some sequences such that~$g^{(k)} \to g$ and~$f_k \to f$ in measure on sets
of finite measure.  By passing to a subsequence we may assume that~$g^{(k)} \to g$ and~$f_k \to f$ almost everywhere, so~$g \in
\mathcal X$ with~$\|g\|_{\mathcal X} \leqslant C$ by the Fatou property.
Now
\begin {equation}
\label {e:sinfdec}
f = (f - f_k) + \mathcal S_\infty g^{(k)} =  (f - f_k) + \mathcal
S_\infty g + \mathcal S_n (g^{(k)} - g) + (\mathcal S_\infty -
\mathcal S_n) \left(g^{(k)} - g\right).
\end {equation}
Let~$E \subset \mathcal M$ be a measurable set of finite measure
and~$\varepsilon > 0$.
The first and the third terms on the right-hand side of~\eqref
{e:sinfdec} converge to~$0$ almost everywhere in~$k$ for any~$n$.
With the help of
the Egoroff theorem and the diagonal process we may choose an increasing
sequence~$n \mapsto k_n$ such that~$\left|\mathcal S_n \left(g^{(k_n)} -
g\right)\right| \leqslant 2^{-n}$ and~$|f - f_{k_n}| \leqslant 2^{-n}$ on a
set~$F \subset E$ such that~$\mu (E \setminus F) < \varepsilon$, and in particular the first and the third terms in~\eqref {e:sinfdec} converge to~$0$ almost
everywhere with~$k = k_n$.
 By the uniform $R$-boundedness
$$
\left\|\left(\mathcal
S_\infty - \mathcal S_n\right) \left(g^{(k_n)} - g\right)\right\|_R \leqslant
\|\mathcal S_\infty - \mathcal S_n\|_{\mathcal X \to R} \left(\left\|g^{(k_n)}\right\|_{\mathcal X} +
\|g\|_{\mathcal X}\right) \to 0,
$$
so the fourth term in~\eqref {e:sinfdec} with~$k = k_n$ converges to~$0$ in
measure on sets of finite measure.  Thus~$f = S_\infty g$ almost everywhere
on~$F$.
Since~$\varepsilon > 0$ and~$E$ are arbitrary, it is easy to see that~$f = S_\infty g$ almost
everywhere, and therefore~$g \in D_{\mathcal X, C} (f)$, which shows that the
graph of~$D_{\mathcal X, C}$ is closed with respect to the convergence in
measure on sets of finite measure.

Now, with~$C > \|f\|_{J (\mathcal X)}$ sets~$D_{\mathcal X, C} (f)$ are
evidently nonempty, and
$$
D_{\mathcal X, \|f\|_{J (\mathcal X)}} (f) = \bigcap_{C > \|f\|_{J (\mathcal
X)}} D_{\mathcal X, C} (f)
$$
is nonempty as an intersection of a centered family of nonempty convex sets that are bounded in
the lattice~$\mathcal X$ with the Fatou property and closed
with respect to the convergence in measure on sets of finite measure by~\cite
[Chapter 10, \S 5, Theorem~3] {kantorovichold}.

\section {The topology of uniform convergence on compact sets}

\label {s:tuccs}

Our methods for establishing certain
properties of interest such as Theorem~\ref {t:jxstability} and its
corollaries are based on a fixed point theorem, and they rely on the closedness
of certain maps in suitable topologies that also make certain bounded sets of
lattices compact.  At present, it is not clear whether it is possible to carry
out these arguments for general spaces~$\Omega$.  Fortunately, at least
for discrete spaces~$\Omega$ there is a natural topology of
uniform convergence on compact sets in~$\mathbb D \times \Omega$ that allows us
to verify the required properties of the maps without much trouble.
\begin {proposition}
\label {latballic}
Let~$X$ be an $r$-convex quasi-normed lattice of measurable functions
on~$\mathbb T \times \Omega$ with a discrete space~$\Omega$ and some~$r > 0$.
Suppose that~$X$ satisfies the Fatou property and property~$(*)$.
Then the closed unit ball~$B_{X_A}$ of~$X_A$ is compact in the topology~$\tau$
of the uniform convergence on all compact sets of~$\mathbb D \times \Omega$.
\end {proposition}
Since~$\tau$ is metrizable, it suffices to verify that for any
sequence~$f_n \in B_{X_A}$ there is a subsequence converging to some~$f \in
B_{X_A}$ in~$\tau$.  We will first prove the claim for Banach lattices~$X$.
Observe that there exists some~$g \in X'$ such that $\|g\|_{X'}  = 1$
and~$g > 0$ a.~e. (see, e.~g., \cite [Proposition~9] {rutsky2010en}).
Lattice~$X'$ also satisfies property~$(*)$ (see \cite [Lemma~2]
{kisliakov2002en}), so there exists some~$\weightw \in X'$ such that~$\weightw
\in X'$, $\weightw > g > 0$ a.~e. and~$\log \weightw (\cdot, \omega) \in
\lclassg {1}$ for almost all~$\omega \in \Omega$.
We may assume that~$\|\weightw\|_{X'} = 1$.
Thus we can construct an outer function
$W = \exp \left(\log \weightw + i \mathcal H [\log \weightw]\right)$ such that
$|W| = \weightw$ almost everywhere.

Let~$\Omega_N \subset \Omega$, $N \in \mathbb N$ be an increasing sequence of
finite sets such that~$\bigcup_N \Omega_N = \Omega$.
We inductively construct a series of increasing sequences
$k \mapsto s_{N, k}$
starting with~$s_{0, k} = k$
such that~$s_{N, k}$ is a subsequence of~$s_{N - 1, k}$.
Sequence~$f_{s_{N - 1, k}} W \chi_{\mathbb T \times \Omega_N}$ belongs to the
unit ball of the space~$\hclass {1} {\mathbb T \times \Omega_N}$, which
is dual to~$\mathrm C \left(\mathbb T \times \Omega_N\right) \slash \mathrm C_A
\left(\mathbb T \times \Omega_N\right)$, so there exists an increasing
subsequence~$s_{N, k}$ of~$s_{N - 1, k}$ such that sequence~$f_{s_{N, k}} W
\chi_{\mathbb T \times \Omega_N}$ converges in the $*$-weak topology to some
$h_N \in \hclass {1} {\mathbb T \times \Omega_N}$,
and therefore~$f_{s_{N, k}} W \chi_{\mathbb D \times \Omega_N} \to h_N$
in~$\tau$.  Functions~$h_M$ and~$h_N$ coincide on~$\mathbb D \times \Omega_M$
for all~$M \leqslant N$, and we may define a function~$h$ by~$h (z, \omega) =
h_N (\omega)$ for~$z \in \mathbb D$ and some~$N$ such that~$\omega \in
\Omega_N$.  For the diagonal sequence~$n' : N \mapsto s_{N, N}$ we have~$f_{n'}
W \to h$ in~$\tau$.

Finally, we need to verify that
$f = W^{-1} h \in B_{X_A}$.
Indeed, by a well-known corollary to the Fatou property
(see, e.~g., \cite [Proposition~10] {rutsky2010en} or~\cite [Proposition~3.3] {rutsky2011en})
there exists a sequence~$\varphi_j$ of finite convex combinations of
$\{f_{n'}\}_{n' > j}$ such that~$\varphi_j \to \varphi$ almost everywhere
on~$\mathbb T \times \Omega$ for some~$\varphi \in B_{X_A}$, and we also
have~$\varphi_j \to f$ in~$\tau$.
Since~$\left\|\varphi_j W (\cdot, \omega)\right\|_{\hclassg {1}} \leqslant 1$
for all~$\omega \in \Omega$, sequence~$\varphi_j W$ satisfies the assumptions of
the Khinchin-Ostrovski\u{\i} theorem (\cite [Chapter~2, \S 8.3] {privaloven}), and
it follows that the boundary values of~$W f$ coincide with~$\varphi W$, thus
the boundary values of~$f$ belong to~$B_X$.

Now suppose that~$X$ is $r$-convex with some $r \geqslant \frac 1 N$ with some
integer $N \geqslant 2$ and~$f_n \in B_{X_A}$.
$X^{\frac 1 N}$ is a Banach lattice because it is $1$-convex and satisfies the
Fatou property, so the conclusion of Proposition~\ref {latballic} applies to it. 
By a similar construction to the above there exists some
$$
\weightw \in \left(X^{\frac 1 N}\right)', \quad \|\weightw\|_{\left(X^{\frac
1 N}\right)'} = 1
$$
and an outer function~$W$ such that $|W| = \weightw$ almost everywhere.
We may raise the inclusion\footnote{To avoid confusion, we remind the reader
that by the definition of a weighted lattice $\lclass {1} {\weightw^{-1}}$ denotes
the ``classical'' weighed space with weight~$\weightw$, i.~e. a lattice with
the norm~$\|h\|_{\lclass {1} {\weightw^{-1}}} = \int |h| \weightw$. }~$X^{\frac
1 N} \subset \lclass {1} {\weightw^{-1}}$ to the power~$N$.
It follows that functions~$F_n = f_n W^N$ belong to the unit
ball of~$\left([\lclassg {1}]^N\right)_A = \hclassg {\frac 1 N}$, and hence they
admit inner-outer factorization~$F_n = I_n G_n$.  Observe that both~$H_{0, n} =
G_n^{\frac 1 N} W^{-1}$ and~$H_{1, n} = I_n H_{0, n}$ belong
to the closed unit ball of~$\left(X^{\frac 1 N}\right)_A$, and~$f_n = H_{0,
n}^{N - 1} H_{1, n}$.  For some subsequence~$n'$ we have $H_{0, n'} \to
h_0$ and~$H_{1, n'} \to h_1$ in~$\tau$ with some~$h_0, h_1$ in the closed unit
ball of~$\left(X^{\frac 1 N}\right)_A$, and $f = h_0^{N - 1} h_1 \in B_{X_A}$ as
claimed.

It is
interesting to note that Proposition~\ref {latballic} allows us to generalize
the equivalence of $\AK$-stability and strong $\AK$-stability to quasi-normed
lattices with discrete space~$\Omega$, although we do not use this
generalization in the present work.  For the case of Banach lattices and
arbitrary~$\Omega$ this was already established in~\cite [Lemma~3] {kisliakov2003}.
\begin {proposition}
\label {saksteq}
Let~$(X, Y)$ be a couple of quasi-normed lattices of measurable
functions on~$\mathbb T \times \Omega$
satisfying the Fatou property and
property~$(*)$.
Suppose also that~$X$ and~$Y$ are $r$-convex with
some~$r > 0$ and~$\Omega$ is discrete.  Then couple~$(X, Y)$ is $\AK$-stable if
and only if it is strongly $\AK$-stable.
\end{proposition}
Indeed, suppose that~$H \in (X + Y)_A$ and~$H = f + g$ with some~$f \in X$
and~$g \in Y$.
Following the proof of~\cite [Lemma~3] {kisliakov2003} we construct
a sequence of outer functions~$\varphi_n$
such that~$|\varphi_n| \leqslant 1$, $\varphi_n \to 1$ in~$\tau$ and~$\varphi_n
H \in X_A + Y_A$.
By the $\AK$-stability there exist some~$F_n \in X_A$ and~$G_n \in Y_A$ such
that~$\varphi_n H = F_n + G_n$,
$\|F_n\|_X \leqslant C \|f\|_X$ and $\|G_n\|_Y \leqslant C \|g\|_Y$.
By Proposition~\ref {latballic} there exists a subsequence~$n'$ such
that~$F_{n'} \to F$ and~$G_{n'} \to G$ in~$\tau$ with some~$F \in X_A$, $G \in
Y_A$, $\|F\|_X \leqslant C \|f\|_X$ and $\|G\|_Y \leqslant C \|g\|_Y$.
But then we also have~$H = F + G$.

\section {Proof of Theorem~\ref {t:jxstability}}

\label {s:jxsp}

We begin by stating a very general fixed point theorem from~\cite {park1998};
for a good general reference on the fixed point theory see, e.~g., \cite {granasdugundji2003}.
This result is not difficult to get hold of for our
purposes, even though the relevant theory is rather complicated
and its key elements may not be familiar or even readily
explainable to an interested reader coming from analysis.

Suppose that $X$ and $Y$ are topological spaces.  A set-valued map 
$T : X \to 2^Y$ is called \emph {closed} if its graph is closed in $X \times Y$.
$T$ is said to be \emph {upper semicontinuous} if for any
closed set $B \subset Y$
its preimage
$$
T^{-1} (B) = \left\{x \in X \mid T (x) \cap B \neq \emptyset \right\}
$$
is also closed.
There is a more natural equivalent definition: $T$ is upper semicontinuous if and only if
for any open set $U \subset Y$ the set
$$
\{x \in X \mid T (x) \subset U\}
$$
is also open.
Thus a composition of of upper semicontinuous maps is also upper semicontinuous.
It is easy to see that if $Y$ is a regular topological
space\footnote {That is, we can separate a point from a closed set not containing it by a couple of open neighbourhoods;
it is well known that any Hausdorff topological vector space is regular.}
and the values of $T$ are closed then $T$ is upper semicontinuous if and only if $T$ is a closed map.
$T$ is called \emph {compact} if the closure of its image $\overline {T (X)}$ is compact in $Y$.
Observe that a composition of compact maps (and even a composition of a compact map with any map)
defined on a Hausdorff compact set is also compact.
For the notion and the definition of an \emph {acyclic} topological space we (by necessity)
refer the reader to \cite {granasdugundji2003}
and to the various algebraic topology textbooks;
in the present work we will only use the simple fact
that convex sets of a topological vector space are acyclic.
$T$ is called an \emph {acyclic map} if $T$ is upper semicontinuous and its values are compact and acyclic.

A nonempty set $X \subset E$ in a linear topological space $E$ is called \emph {admissible} (in the sense of Klee)
if for any compact set $K \subset X$ and any open set $V \subset E$, $0 \in V$, there exists
a continuous map $h : K \to X$ such that
$x - h (x) \in V$ for all $x \in K$ and $h (K)$ is contained in a finite-dimensional subspace $L \subset E$.
In other words, $X$ is admissible if any compact set $K \subset X$ can be continuously and uniformly approximated
by a family of finite-dimensional sets of $X$.
In particular, any non\-empty convex set of a locally convex linear topological space is admissible.

Let $X$ be a nonempty convex set in a linear topological space $E$,
and let $Y$ be another linear topological space.
A set $P \subset X$ is called a \emph {polytope} if $P$ is the convex hull of a finite set in $X$.
A map $F : X \to 2^Y$ belongs to the \emph {``better'' admissible class} $\mathfrak B (X, Y)$ if and only if
for any polytope $P \subset X$ and any continuous function $f : F (P) \to P$ the composition
$f \circ F|_P : P \to P$ has a fixed point.
Observe that \emph {admissibility} refers in this notion to the existence of fixed points in a restricted sense.
The class $\mathfrak B (X, Y)$ encompasses a large number of particular classes
of maps that are known to have fixed points.
In the present work we will only
use the fact that this class contains finite compositions of acyclic maps.
The corresponding fixed point theorem was established in~\cite {powers1972} (see also \cite [\S 19.9] {granasdugundji2003}
for the statement in context), and 
it is possible to use it directly with minor adaptations, similarly to how
the Fan--Kakutani fixed point theorem is derived from the Kakutani fixed point
theorem.
The result~\cite {park1998}, however, allows us to keep the necessary topological explanations to a minimum.
\begin {theorem} [{\cite [Corollary~1.1] {park1998}}]
\label {parkfptc}
Let $E$ be a Hausdorff topological vector space, and let $X \subset E$ be an admissible convex set.
Then any closed compact map $\Phi \in \mathfrak B (X, X)$ has a fixed point.
\end {theorem}

We will carry out some arguments based on a fixed point theorem applied to maps
acting on sets of some majorants in lattices.  In many cases it suffices to
endow the sets of logarithms of such majorants with the weak topology of a weighted space~$\lclass {2}
{\omega}$; see, e.~g., \cite [Proposition~30] {rutsky2016en}.  However,
it is not clear if the convergence in this topology
implies the convergence in~$\tau$, since in general~$\log \omega \notin
\lclassg {1}$ in the first variable.  Otherwise property~$(*)$ could have been
improved to majorants with logarithms in~$\lclassg {2}$, but it is easy to find
Orlicz spaces with property~$(*)$ that do not admit majorants with logarithms
in~$\lclassg {p}$ for any~$p > 1$.  Fortunately, the sets of majorants still
turn out to be compact in the weak topology of~$\lclassg {1}$.  Some care needs
to be taken, however, since it is not clear if these sets are
separable in this topology.

\begin {lemma}
\label {logl1isl2}
Let $f \in \lclassg {1}$ and $f \geqslant 1$ almost everywhere.
Then $\log f \in \lclassg {2}$ and  $\|\log f\|_{\lclassg {2}} \leqslant 2 \|f\|_{\lclassg {1}}^{\frac 1 2}$.
\end {lemma}
We only need to observe that $\log \left(f^{\frac 1 2}\right) \leqslant f^{\frac
1 2}$, and so $$
\int (\log f)^2 = 4 \int \left(\log f^{\frac 1 2}\right)^2 \leqslant 4 \int f.
$$

\begin {proposition}
\label {luscompact}
Let~$X$ be a Banach lattice of measurable functions on $\mathbb T \times \Omega$
with discrete~$\Omega$, and let~$f \in X$.
Suppose that~$X$ satisfies the Fatou property and property~$(*)$.  Suppose also
that $\log f (\cdot, \omega) \in \lclassg {1}$ for all~$\omega \in \Omega$.
Then for all~$A > 0$ sets
\begin {equation*}
V_{X, f, A} = \{\log g \in \lclassg {1} \mid g \geqslant f, \|g\|_X \leqslant A\}
\end {equation*}
are compact in the topology of weak convergence in~$\lclass {1} {\mathbb T
\times \{\omega\}}$ for all~$\omega \in \Omega$.
\end {proposition}
First, suppose that~$\Omega$ is trivial, so we only have one variable.
It is well known that the positive part of the unit ball of a Banach lattice is logarithmically convex,
so~$V_{X, f, A}$ is a convex set.
Observe first that by the Fatou property~$V_{X, f, A}$ is closed with respect to
the convergence in measure on sets of finite measure, thus~$V_{X, f, A}$ is also
closed in~$\lclassg {1}$ and therefore weakly closed.
Now, by the Dunford-Pettis theorem it suffices to prove that the set $V_{X, f, A}$
is bounded and uniformly absolutely continuous.
Let $B$ be a measurable set of $\mathbb T$.
We take some~$\weightw \in X'$ as in the proof of Proposition~\ref {latballic}.
Then $\{\weightw g \mid \log g \in V_{X, f, A}\}$ is a bounded set in~$\lclassg
{1}$, and by Lemma~\ref {logl1isl2} we have
$\left\|\log^+ [\weightw g]\right\|_{\lclassg {2}} \leqslant 2 A^{\frac 1 2}$
for all $\log g \in V_{X, f, A}$.
Thus
$$
\int_B \log^+ [\weightw g] \leqslant \left\|\log^+ [\weightw
g]\right\|_{\lclassg {2}} |B|^{\frac 1 2} \to 0 $$
as $|B| \to 0$, and this convergence is uniform in~$\log g \in V_{X, f, A}$.
On the other hand, we also have
$$
\int_B \log^- [\weightw g] \geqslant \int_B \log^- [\weightw f] \to 0
$$
as~$|B| \to 0$ uniformly in~$\log g \in V_{X, f, A}$.
Therefore,
$$
\left\{\log [\weightw g] \mid \log g \in V_{X, f, A}\right\} = \log \weightw +
V_{X, f, A} $$
is bounded and uniformly absolutely continuous, which implies its relative
com\-pac\-tness in the weak topology of $\lclassg {1}$.
It follows that $V_{X, f, A}$ is also compact.

Now, for arbitrary discrete~$\Omega$ we consider the restricted
lattices~$X_\omega = \{h (\cdot, \omega) \mid h \in X\}$ with the corresponding
norm~$\|g\|_{X_\omega} = \|g\|_{X}$, $\omega \in \Omega$.  They
also satisfy the assumptions of Proposition~\ref {luscompact}. By the Tychonoff
theorem space~$V = \prod_{\omega \in \Omega} V_{X_\omega, f (\cdot, \omega), A}$
is compact with the product topology.
It suffices to show that $V_{X, f, A} \subset V$ is closed in~$V$.
Observe that the map
\begin {equation}
\label {phi1def}
\Phi_0 (\log g) (z, \omega) =
\exp \left(\frac 1 {2 \pi} \int_0^{2 \pi} \frac {e^{i \theta} + z} {e^{i \theta}
- z} \log g \left(e^{i \theta}, \omega\right) d\theta \right)
\end {equation}
defined for all $\log g (\cdot, \omega) \in \lclassg {1}$, $z \in \mathbb D$ and
$\omega \in \Omega$ is continuous as a map from~$V_{X, f, A}$ to~$A B_{X_A}$
with topology~$\tau$.  The integral in~\eqref {phi1def} is the convolution with
the Schwarz kernel, and~$\Phi_0 (\log g)$ is an outer function with the boundary
values satisfying $|\Phi_0 (\log g)| = g$ almost everywhere.
If~$\log g_\alpha \in V_{X, f, A}$ is a net converging to some~$\log g \in V$
in~$V$ then~$\Phi_0 (\log g_\alpha) \to \Phi_0 (\log g)$,
and by Proposition~\ref {latballic} $\Phi_0 (\log g) \in A B_{X_A}$,
so~$\|g\|_X \leqslant A$.  On the other hand, since~$|\Phi_0 (\log g_\alpha)|
\geqslant |\Phi_0 (\log f)|$ on~$\mathbb T$, we also have~$|\Phi_0 (\log
g_\alpha) (z, \omega)| \geqslant |\Phi_0 (\log f) (z, \omega)|$ for all~$z \in
\mathbb D$ and~$\omega \in \Omega$. Passing to the limit in~$\alpha$
yields $|\Phi_0 (\log g) (z, \omega)| \geqslant |\Phi_0 (\log f) (z,
\omega)|$, and passing then to the boundary values shows that~$g \geqslant f$.
Therefore, $\log g \in V_{X, f, A}$, and~$V_{X, f, A}$ is indeed closed.  

We now begin the proof of Theorem~\ref {t:jxstability}.  Suppose that under its
assumptions~$f_0 \in \left[J (\mathcal X)\right]_A$ with norm~$1$, and~$J
(\mathcal X)$ is $\nplus$-stable with constant~$C$.
By property~$(*)$ there exists some~$f \geqslant |f_0|$, $\|f\|_{J (\mathcal X)}
\leqslant 2$ and~$\log f (\cdot, \omega) \in \lclassg {1}$ for all~$\omega$.
We define a set-valued map~$\Phi_1 : 4 B_{[J (\mathcal X)]_A} \to 2^{4 C
B_{\mathcal X_A}}$ by $$
\Phi_1 (g) = \left\{ \{g_j\}_{j \in I} \in \mathcal X_A \mid g = \sum_{j \in I}
g_j, \left\|\{g_j\}_{j \in I}\right\|_{\mathcal X} \leqslant 4 C\right\}, \quad g \in
4 B_{[J (\mathcal X)]_A}.
$$
It has nonempty convex values.  Let us show that~$\Phi_1$ is upper
semicontinuous with respect to~$\tau$.  Indeed, suppose that~$g_k \in 4 B_{X_A}$ and~$h_k \in
\Phi_1 (g_k)$ are such that~$g_k \to g$ and~$h_k \to h$ in~$\tau$.
By~\cite [Proposition~10] {rutsky2010en} and the convexity of the graph
of~$\Phi_2$ we may replace them with a sequence of convex combinations such that
additionally~$g_k \to g$ and~$h_k \to h$ almost everywhere.  By Proposition~\ref
{p:jxsumcm} it follows that~$h \in \Phi_1 (g)$.

We also define a set-valued map~$\Phi_2 : 4 C B_{\mathcal X_A} \to
2^{V_{J (\mathcal X), f, 4}}$ by
$$
\Phi_2 (\{h_j\}_{j \in I}) = \left\{ \log \weightw \in V_{J (\mathcal X), f, 4}
\mid \weightw \geqslant f + \frac 1 {2 C} \sum_{j \in I} |h_j|\right\}
$$
for $\{h_j\}_{j \in I} \in 4 C B_{\mathcal X_A}$.  This map also has nonempty
convex values.  To verify its upper semicontinuity,
suppose that~$h^{(k)} = \left\{h_j^{(k)}\right\}_{j
\in I} \in 4 C B_{\mathcal X_A}$ and~$\log \weightw_k \in \Phi_2
\left(\{h_j\}_{j \in I}\right)$ are such that~$h^k \to h$ in~$\tau$ and~$\log
\weightw_k (\cdot, \omega) \to \log \weightw (\cdot, \omega)$ in the weak
topology of~$\lclassg {1}$ for all~$\omega \in \Omega$. Let~$n \in
\mathbb N$.  Estimate
\begin {equation}
\label {e:shest}
\log \left(|\Phi_0
(\log f)| + \frac 1 {2 C} \sum_{j \in I \cap [-n, n]}
\left|h_j^{(k)}\right|\right) \leqslant \log |\Phi_0 (\log \weightw_k)|
\end {equation}
is satisfied on~$\mathbb T \times \Omega$.  By~\cite [Proposition~2.2]
{graczyketal2010} the function on the left-hand side of~\eqref {e:shest} is
subharmonic, and the function on the right-hand side of~\eqref {e:shest} is
harmonic, so~\eqref {e:shest} is also satisfied on~$\mathbb D \times \Omega$.
Passing to the limit in~$k$, passing to the boundary values and then passing to
the limit in~$n$ shows that~$\log \weightw \in \Phi_2 (h)$.

Now we define the composition map~$\Phi = \Phi_2 \circ \Phi_1 \circ \Phi_0$,
which belongs to the class $\mathfrak B \left(V_{J (\mathcal X), f, 4}, V_{J (\mathcal X), f, 4}\right)$
as a finite composition of maps taking acyclic values.
$\Phi$ is closed and compact as a composition of compact upper
semicontinuous maps.
Thus by Theorem~\ref {parkfptc} there exists some~$\log \weightw \in V_{J
(\mathcal X), f, 4}$ such that $\log \weightw \in \Phi (\log \weightw)$.
This also implies that there exist some~$\{h_j\}_{j \in I} \in 4 C B_{\mathcal
X_A}$ satisfying~$\Phi_0 (\log \weightw) = \sum_{j \in I} h_j$ and~$\sum_{j \in
I} |h_j| \leqslant 2 C |\Phi_0 (\log \weightw)|$.  Therefore,
functions~$\varphi_j = \frac {h_j} {\Phi_0 (\log \weightw)}$, $j \in I$ satisfy
$\left\|\{\varphi_j\}_{j \in I}\right\|_{\hclass {\infty} {\lsclass {1}}}
\leqslant 2 C$, $\sum_{j \in I} \varphi_j = 1$, and~$|f_0 \varphi_j| \leqslant
\frac {|f_0|} {\weightw} |h_j| \leqslant |h_j|$ almost everywhere on~$\mathbb T
\times \Omega$ for all~$j \in I$, which shows that~$\{\varphi_j\}_{j \in I}$
provide the $\nplus$-stability for~$f_0$ with constant~$4 C$ in the sense of
Definition~\ref {d:jxstability}.  We note in passing that the constant can be
improved to~$2 C$ if one takes advantage of the fact that the constant in
property~$(*)$ can be made arbitrarily close to~$1$ (see~\cite [Lemma~2.2]
{kalton1994}).

\section {Sufficiency of condition~\ref {t:c:xyincl} for AK-stability}

\label {s:socxy}

The following result is in a certain way a natural development of Corollary~\ref
{akbakeq}.  Since~$X^{1- \theta} Y^\theta$ is a space of type~$\mathcal
C_\theta (X, Y)$, it also provides a nice
direct and self-contained link between $\nplus$-stability with respect
to~$(\cdot, \cdot)_{\theta, \infty}$ and $\AK$-stability.  Observe that the
$\AK$-stability of~$(X, Y)$ naturally implies the $\nplus$-stability of this
couple with respect to~$(\cdot, \cdot)_{\theta, \infty}$, which in turn implies
the inclusion~\eqref {e:xyincl} below.
Moreover, after a more direct proof based on the Powers fixed point theorem we
will show how a natural modification of this technique (building on some of the ideas
from~\cite {rutsky2017}) allows one to prove Corollary~\ref {akbakeq} using only the Fan--Kakutani fixed point theorem.
\begin {proposition}
\label {p:xyincl}
Suppose that~$(X, Y)$ is a couple of Banach lattices of
measurable functions on~$\mathbb T \times \Omega$ with a discrete space~$\Omega$
satisfying the Fatou property and property~$(*)$.
Then
\begin {equation}
\label {e:xyincl}
\left(X^{1 - \theta} Y^\theta\right)_A \subset \left(X_A,
Y_A\right)_{\theta, \infty}
\end{equation}
with some~$0 < \theta < 1$ if and only if
couple~$(X, Y)$ is boundedly $\AK$-stable.
\end {proposition}

Indeed, suppose that~\eqref {e:xyincl} is true and we are given some~$f \in X$
and~$g \in Y$.  The argument that follows is similar to the proof of
Theorem~\ref {t:jxstability}, although the details (that we repeat for clarity)
are somewhat simpler.
We may
assume that~$f$ and~$g$ satisfy~$\log f (\cdot, \omega), \log g (\cdot, \omega) \in \lclassg {1}$ for all~$\omega \in \Omega$.  Let~$D_X = V_{X, f, 2 \|f\|_X}$ and~$D_Y = V_{Y, g, 2
\|g\|_Y}$ be the sets defined in Proposition~\ref {luscompact}. 

Let~$\log u \in D_X$ and~$\log v \in D_Y$.  Then~$w = u^{1 - \theta} v^\theta
\in X^{1 - \theta} Y^{\theta}$ with norm at most~$2 \|f\|_X^{1 - \theta}
\|g\|_Y^\theta$, and we construct the outer function~$W = \Phi_0 (\log w)$,
$|W| = w$ on~$\mathbb T \times \Omega$, where~$\Phi_0$ is the map defined by~\eqref {phi1def}.
Thus~$W \in \left(X^{1 - \theta} Y^{\theta}\right)_A$ with the same estimate for
the norm as~$w$, and from~\eqref {e:xyincl} it follows that~$W = F + G$
with some~$F \in X_A$ and~$G \in Y_A$
satisfying~$\|F\|_X \leqslant C t^\theta {\|f\|\strut}_X^{1 - \theta}
{\|g\|\strut}_Y^\theta$ and~$\|G\|_Y \leqslant C t^{\theta - 1}
{\|f\|\strut}_X^{1 - \theta} {\|g\|\strut}_Y^\theta$ for all~$t > 0$ with
some~$C$ independent of~$f$, $g$ and~$t$.
Choosing~$t = \frac {\|f\|_X} {\|g\|_Y}$ yields~$\|F\|_X \leqslant C \|f\|_X$
and~$\|G\|_Y \leqslant C \|g\|_Y$.  Let~$C_X = C \|f\|_X$ and~$C_Y = C \|g\|_Y$.
We see that a set-valued map~$\Phi_1 : D_X \times D_Y \to 2^{C_X B_{X_A}
\times C_Y B_{Y_A}}$ defined by
\begin {multline*}
\Phi_1 (\log u, \log v) =
\\
\left\{ (F, G) \mid F \in C_X B_{X_A}, G \in C_Y B_{Y_A},
\Phi_0 \left(\log \left[u^{1 - \theta} v^\theta\right]\right) = F + G \right\}
\end {multline*}
takes nonempty convex values for all~$\log u \in D_X$, $\log v \in D_Y$.
We also define a set-valued map~$\Phi_2 : C_X B_{X_A} \times C_Y B_{Y_A} \to
2^{D_X \times D_Y}$ by
\begin {multline*}
\Phi_2 (F, G) = \left\{ (\log u_1, \log v_1) \mid \|u_1\|_X \leqslant 2 \|f\|_X,
\|v_1\|_Y \leqslant 2 \|g\|_Y, \right.
\\
\left.
u_1 \geqslant f \vee \frac 1 C |F|, v_1 \geqslant g \vee \frac 1 C |G|
\right\}, \quad F \in C_X B_{X_A}, G \in C_Y B_{Y_A},
\end {multline*}
that also takes nonempty convex values.

We endow~$D_X$ and~$D_Y$ with the topology of weak convergence in~$\lclass
{1} {\mathbb T \times \{\omega\}}$ for all~$\omega \in \Omega$ and~$C B_{X_A}$,
$C B_{Y_A}$ with the topology~$\tau$ of uniform convergence on compact sets
of~$\mathbb D \times \Omega$, which turns them into compact convex sets in the
respective locally convex linear topological spaces (see
Proposition~\ref {luscompact} and Proposition~\ref {latballic}).

It is easy to see that~$\Phi_1$ is upper semicontinuous.
Let us show that~$\Phi_2$ is also upper semicontinuous.
Suppose that~$F_k \in C B_{X_A}$, $G_k \in C B_{Y_A}$ and~$(\log u_k, \log v_k)
\in \Phi_2 (F_k, G_k)$ are such that~$F_k \to F$, $G_k \to G$, $\log u_k \to
\log u$ and~$\log v_k \to \log v$ in the respective spaces.  We construct outer
functions~$U_k = \Phi_0 (\log u_k)$, $V_k = \Phi_0 (\log v_k)$, $U = \Phi_0
(\log u)$, $V = \Phi_0 (\log v)$, $\varphi = \Phi_0 (\log f)$ and~$\psi = \Phi_0
(\log g)$.  Then~$|U_k| \geqslant |\varphi| \vee \frac 1 C |F_k|$ and~$|V_k|
\geqslant |\psi| \vee \frac 1 C |G_k|$ on~$\mathbb D \times \Omega$.  Passing to
the limit in~$k$ yields~$|U| \geqslant |\varphi| \vee \frac 1 C |F|$ and~$|V|
\geqslant |\psi| \vee \frac 1 C |G|$ on~$\mathbb D \times \Omega$, and therefore
also almost everywhere on~$\mathbb T \times \Omega$, so~$u \geqslant f \vee
\frac 1 C |F|$, $v \geqslant g \vee \frac 1 C |G|$ and~$(\log u, \log v) \in
\Phi_2 (F, G)$.

Now we define the composition map~$\Phi = \Phi_2 \circ \Phi_1$.  It belongs
to~$\mathfrak B (D_X \times D_Y, D_X \times D_Y)$ as a finite composition of
maps taking acyclic values, and it is compact as a finite composition of compact
maps.  Therefore, by Theorem~\ref {parkfptc} there exist some~$\log u \in D_X$,
$\log v \in D_Y$ such that~$(\log u, \log v) \in \Phi (\log u, \log v)$.
This means that
$W = \Phi_0 \left(\log \left[u^{1 - \theta} v^\theta\right]\right) = F + G$
for some~$F \in C_X B_{X_A}$, $G \in C_Y B_{Y_A}$ satisfying~$|F| \leqslant C u$
and~$|G| \leqslant C v$ on~$\mathbb T \times \Omega$.

Let~$V = \frac F W$ and~$U = \frac G W$.  Then~$V + U = 1$, $\chi_{\{u \leqslant
v\}} |V| \leqslant \chi_{\{u \leqslant v\}} \frac {|F|} {u} \leqslant C$, and,
similarly, $\chi_{\{u > v\}} |U| \leqslant C$ almost everywhere on~$\mathbb T
\times \Omega$.  Therefore, also~$\chi_{\{u > v\}} |V| = \chi_{\{u > v\}} |1 -
U| \leqslant C + 1$, and thus~$V \in \hclassg {\infty}$ with norm at most~$C +
1$.  Observe that~$|V| g^\theta = \frac {|F| g^\theta} {u^{1 - \theta} v^\theta}
\leqslant C u^\theta$ and, similarly, $|1 - V| f^{1 - \theta} \leqslant C v^{1 -
\theta}$.  This implies that~$V$ satisfies the conditions of
Lemma~\ref {bakstabilitygen} with $\alpha = \frac 1 \theta$ and~$\beta = \frac 1
{1 - \theta}$, so~$(X, Y)$ is indeed boundedly $\AK$-stable.

Now we will show how Proposition~\ref {p:xyincl} can be derived from the
Fan--Kakutani fixed point theorem using a suitable approximation. 
Let~$\Omega_N$ be as in the proof of Proposition~\ref {latballic}, $0 < r_N <
1$, $\varepsilon_N > 0$, $r_N \to 1$ and~$\varepsilon_N \to 0$.
It is easy to see that a set-valued map~$\tilde \Phi_1^{(N)} : D_X \times D_Y
\to 2^{C_X B_{X_A} \times C_Y B_{Y_A}}$ defined by
\begin {multline*}
\tilde \Phi_1^{(N)} (\log u, \log v) =
\left\{ (F, G) \mid F \in C_X B_{X_A}, G \in C_Y B_{Y_A},\right.
\\
\left.
\left|\Phi_0 \left(\log \left[u^{1 - \theta} v^\theta\right]\right) - (F +
G)\right| < \varepsilon_N \text { on $r_N \overline {\mathbb D} \times
\Omega_N$} \right\}
\end {multline*}
is lower semicontinuous\footnote{In fact, $\Phi_1^{(N)}$ has open graph.}:
if~$\log u_\alpha \in D_X$, $\log v_\alpha \in D_Y$ are some nets converging to
some functions~$\log u$, $\log v$ and~$(F, G) \in \tilde \Phi_1^{(N)} (\log u,
\log v)$ then~$\Phi_0 \left(\log \left[u_\alpha^{1 - \theta}
v_\alpha^\theta\right]\right)$ converges in~$\tau$ to~$\Phi_0 \left(\log \left[u^{1 - \theta}
v^\theta\right]\right)$, so 
$\left|\Phi_0 \left(\log \left[u_\alpha^{1 - \theta}
v_\alpha^\theta\right]\right) - (F + G)\right| < \varepsilon_N$
for all~$\alpha \succ \beta$ with some~$\beta$ on~$r_N \overline {\mathbb D}
\times \Omega_N$.  Therefore, the closure of this map
\begin {multline*}
\overline {\tilde \Phi_1^{(N)}} (\log u, \log v) =
\left\{ (F, G) \mid F \in C_X B_{X_A}, G \in C_Y B_{Y_A},\right.
\\
\left.
\left|\Phi_0 \left(\log \left[u^{1 - \theta} v^\theta\right]\right) - (F +
G)\right| \leqslant \varepsilon_N \text { on $r_N \overline {\mathbb D} \times
\Omega_N$} \right\}
\end {multline*}
is also lower semicontinuous, and it takes nonempty convex compact values.
By the Michael selection theorem \cite {michael1956} there exists a continuous
selection~$\Phi_1^{(N)} : D_X \times D_Y \to C_X B_{X_A} \times C_Y B_{Y_A}$
of~$\overline {\tilde \Phi_1^{(N)}}$, that is, $\Phi_1^{(N)} (\log u, \log v)
\in \overline {\tilde \Phi_1^{(N)}} (\log u, \log v)$ for all~$\log u \in D_X$,
$\log v \in D_Y$.

We now proceed as before with a set-valued map~$\Phi^{(N)} = \Phi_2 \circ
\Phi_1^{(N)}$ in place of~$\Phi$.  This map takes convex values and is upper
semicontinuous.
By the Fan--Kakutani theorem \cite {fanky1952} (which is a much less involved
particular case of Theorem~\ref {parkfptc}) $\Phi^{(N)}$ has some fixed
points~$\log u_N \in D_X$, $\log v_N \in D_Y$, $(\log u_N, \log v_N) \in
\Phi^{(N)} (\log u_N, \log v_N)$.  Let~$(F_N, G_N) = \Phi_1^{(N)} (\log u_N,
\log v_N)$. By the compactness of~$D_X \times D_Y$ and~$C_X B_{X_A} \times C_Y
B_{Y_A}$ these sequences have some limit points~$\log u \in D_X$, $\log v \in
D_Y$, $F \in C_X B_{X_A}$, $G \in C_Y B_{Y_A}$ respectively, and it is easy to
see that~$(F, G) \in \Phi_1 (\log u, \log v)$, so~$(\log u, \log v)$ is a fixed point of the
original map~$\Phi$, and the proof may proceed as before. 

\section {The scale of the bounded AK-stability can be extended}

\label {s:trtbakst}

The following result is a generalization of Proposition~\ref {linfbmorbakst}.
It
is also a key component for both the necessity and the sufficiency of weak-type $\BMO$-regularity for the stability of the real interpolation in
Theorem~\ref {goodrealinterp}.

\begin {theorem}
\label {t:erbakst}
Let~$(X_0, X_1)$ and~$(Y_0, Y_1)$ be two couples of $r$-convex
quasi-normed lattices of measurable functions on~$\mathbb T \times \Omega$ with
some~$r > 0$ and a discrete~$\Omega$ satisfying the Fatou property and
property~$(*)$.
Suppose that~$Y_j$ is of type~$\mathcal C_{\theta_j} (X_0, X_1)$, $j \in \{0, 1\}$ with some~$0 < \theta_0 < \theta_1 < 1$.  If~$(Y_0,
Y_1)$ is boundedly $\AK$-stable then so is~$(X_0, X_1)$.
\end{theorem}

First of all, we may assume that~$r \leqslant 1$.  By Proposition~\ref
{bakrscale}, Proposition~\ref {p:riscale} and the inclusion~$(X, Y)_{\theta_j, 1} \subset (X, Y)_{\theta_j, \frac 1
r}$ we may raise all lattices to the power~$\frac 1 r$ and thus assume that they
are all Banach.

Let~$Z_j = (Y_0, Y_1)_{\delta_j, \infty}$, $j \in \{0, 1\}$
with some~$0 < \delta_0 < \delta_1 < 1$.
Then by the reiteration theorem~$Z_j = (X_0, X_1)_{\alpha_j, \infty}$ with
some~$0 < \alpha_0 < \alpha_1 < 1$, and this couple is $\AK$-stable by~\cite
[Lemma~1.1] {kisliakov1999}.  Then it is boundedly $\AK$-stable with a
constant~$C$ by Corollary~\ref {akbakeq}.

The proof now follows essentially the same idea as the proof of
Proposition~\ref {linfbmorbakst}.  Let~$\alpha_0 < \beta_0 < \beta_1 <
\alpha_1$.  We will first show that couple~$\left(E_0, E_1\right)$ with~$E_j =
X_0^{1 - \beta_j} X_1^{\beta_j}$, $j \in \{0, 1\}$ is boundedly $\AK$-stable.
Suppose that~$f_j \in E_j$, $j \in \{0, 1\}$ are some nonnegative functions
such that~$\|f_j\|_{E_j} = 1$, and let~$t > 0$.  With the help of
property~$(*)$ we may assume that~$\log f_j (\cdot, \omega) \in \lclassg {1}$
for all~$\omega \in \Omega$.

Let~$\gamma_0 = 1 - \alpha_0  - (1 - \beta_0) \frac {\alpha_0}
{\beta_0}$, $\gamma_1 = \alpha_1 - \frac {1 - \alpha_1} {1 - \beta_1} \beta_1$,
$\zeta_0 = 1 - (1 - \alpha_1) \frac {1 - \beta_0} {1 - \beta_1}$
and~$\zeta_1 = \alpha_0 \frac {\beta_1} {\beta_0}$.  Note that the arguments
that follow are symmetric with respect to interchanging~$X_0$ and~$X_1$ and simultaneously
replacing~$\alpha_j$ with~$1 - \alpha_{1 - j}$, $j \in \{0, 1\}$, and
the same is true for~$\beta_j$ and~$\zeta_j$, so it suffices to verify these
arguments and computations for one side only.

It is easy to
see that~$0 < \gamma_0, \gamma_1 < 1$
and~$\zeta_0 = (\beta_0 - \beta_1) \frac {1 - \alpha_1} {1 - \beta_1} + \alpha_1 =
\beta_0 \frac {1 - \alpha_1} {1 - \beta_1} + \gamma_1$,
$\zeta_1 = (\beta_1 - \beta_0) \frac {\alpha_0} {\beta_0} + \alpha_0 = (1 -
\beta_1) \frac {\alpha_0} {\beta_0} + 1 - \gamma_0$,
so in particular~$\alpha_0 < \zeta_0, \zeta_1 < \alpha_1$ and~$\gamma_1
< \zeta_0$, $\gamma_0 < 1 - \zeta_1$.
We take some~$\omega_j \in X_j^{\gamma_j}$ with norm~$1$
such that~$\omega_j > 0$ almost everywhere, $j \in \{0, 1\}$ (see, e.~g., \cite [Proposition~9]
{rutsky2010en}).  By making use of property~$(*)$ we may also assume
that~$\log \omega_j (\cdot, \omega) \in \lclassg {1}$ for all~$\omega \in
\Omega$.
Let~$D_j = V_{X_j^{\gamma_j}, \omega_j, 2}$, $j \in \{0, 1\}$ be the sets
defined in Proposition~\ref {luscompact}.
Observe that for any~$\log u_j \in D_j$ we
have~$g_0 = f_0^{\frac {\alpha_0} {\beta_0}} u_0 \in X_0^{1 - \alpha_0} X_1^{\alpha_0} \subset Z_0$ 
and~$g_1 = f_1^{\frac {1 - \alpha_1} {1 - \beta_1}} u_1 \in X_0^{1 - \alpha_1}
X_1^{\alpha_1} \subset Z_1$.  The norms of these functions in these spaces
are at most~$c$, a constant independent of~$f_0$, $f_1$ and~$t$. 
By the bounded $\AK$-stability of~$(Z_0, Z_1)$ there exists some~$U \in \hclassg {\infty}$ with norm at most~$C$ such that~$\|U g_1\|_{Z_0}
\leqslant c C t$ and~$\|(1 - U) g_0\|_{Z_1} \leqslant c C t^{-1}$.

Let~$F_j = X_0^{1 - \zeta_j} X_1^{\zeta_j}$, $j \in \{0, 1\}$.  Since~$F_j$ is a
space of type~$\mathcal C_{\eta_j} (Z_0, Z_1)$ with~$\eta_j = \frac {\zeta_j - \alpha_0} {\alpha_1 - \alpha_0}$,
we have
$\|U g_1\|_{F_0} \leqslant c_1 \|U g_1\|_{Z_0}^{1
- \eta_0} \|U g_1\|_{Z_1}^{\eta_0} \leqslant c_2 t^{1 - \eta_0}$, and,
similarly,
$\|(1 - U) g_0\|_{F_1} \leqslant c_3 t^{-\eta_1}$
with some~$c_1$, $c_2$ and~$c_3$ independent of~$f_0$, $f_1$ and~$t$.
Thus, a set-valued map~$\Phi_1 : D_0 \times D_1 \to 2^{D_0 \times D_1 \times C
B_{\hclassg {\infty}}}$ defined by
\begin {multline*}
\Phi_1 (\log u_0, \log u_1) = \left\{ (\log u_0, \log u_1, U) \mid  U \in
\hclassg {\infty},  \|U\|_{\hclassg {\infty}} \leqslant C,\right.
\\
\left\|U f_1^{\frac {1 - \alpha_1} {1 - \beta_1}} u_1\right\|_{F_0} \leqslant
c_2 t^{1 - \eta_0}, \left.
\left\|(1 - U) f_0^{\frac {\alpha_0} {\beta_0}} u_0\right\|_{F_1} \leqslant c_3
t^{-\eta_1} \right\}
\end {multline*}
takes nonempty convex values.

Now suppose that~$U \in \Phi_1 (\log u_0, \log u_1)$.  From the definition
of the pointwise lattice products
$$
F_0 = \left(X_0^{1 - \zeta_0} X_1^{\zeta_0 -
\gamma_1}\right) X_1^{\gamma_1},
$$
$$F_1 = \left(X_0^{1 - \zeta_1 - \gamma_0}
X_1^{\zeta_1}\right) X_0^{\gamma_0}
$$
it follows that there exist some nonnegative
functions~$v_j \in B_{X^{\gamma_j}}$, $j \in \{0, 1\}$ such that $$
\left\|U f_1^{\frac {1 - \alpha_1} {1 - \beta_1}} u_1 v_1^{-1}\right\|_{X_0^{1- \zeta_0} X_1^{\zeta_0 - \gamma_1}} \leqslant 2 c_2
t^{1 - \eta_0},
$$
$$
\left\|(1 - U) f_0^{\frac {\alpha_0} {\beta_0}}
u_0 v_0^{-1} \right\|_{X_0^{1- \zeta_1 - \gamma_0} X_1^{\zeta_1}} \leqslant 2
c_3 t^{-\eta_1}.
$$
By replacing~$v_j$ with~$v_j \vee \omega_j$, $j \in \{0, 1\}$ we may assume that~$\log v_j \in D_j$
at the same time as these estimates hold true.
Thus,
a set-valued map~$\Phi_2 : \Phi_1 (D) \to 2^{D_0 \times D_1}$
defined by
\begin {multline*}
\Phi_2 (\log u_0, \log u_1, U) = \left\{
(\log v_0, \log v_1) \mid \right.
\\
\left\|U f_1^{\frac {1 - \alpha_1} {1 - \beta_1}} u_1 v_1^{-1}\right\|_{X_0^{1- \zeta_0} X_1^{\zeta_0 - \gamma_1}} \leqslant 2 c_2
t^{1 - \eta_0},
\\
\left.
\left\|(1 - U) f_0^{\frac {\alpha_0} {\beta_0}}
u_0 v_0^{-1} \right\|_{X_0^{1- \zeta_1 - \gamma_0} X_1^{\zeta_1}} \leqslant 2
c_3 t^{-\eta_1}  
\right\}
\end {multline*}
takes nonempty values for~$(\log u_0, \log u_1, U) \in D$ that are convex.

Now let~$\Phi = \Phi_2 \circ  \Phi_1$.
Observe that~$\Phi \in \mathfrak B \left(D_0 \times D_1, D_0 \times D_1\right)$
as a finite composition of maps taking acyclic values.
We will now show that this map has a fixed
point by verifying that it satisfies the assumptions of Theorem~\ref {parkfptc}.

We endow~$D_0$ and~$D_1$ with the topology of weak convergence in~$\lclass {1}
{\mathbb T \times \{\omega\}}$ for all~$\omega \in \Omega$ and~$C
B_{\hclassg {\infty}}$ with the topology~$\tau$ of uniform convergence on compact sets
in~$\mathbb D \times \Omega$, which turns them into compact convex sets in
the respective locally convex linear topological spaces (see
Proposition~\ref {luscompact} and Proposition~\ref {latballic}).

Let us verify that map~$\Phi_1$ is upper semicontinuous.  Suppose
that we are given some nets~$\log u_j^{(\alpha)} \in D_j$, $\left(\log
u_0^{(\alpha)}, \log u_1^{(\alpha)}, U_\alpha
\right) \in \Phi_1 \left(\log u_0^{(\alpha)}, \log u_1^{(\alpha)}\right)$,
and~$\log u_j^{(\alpha)} \to \log u_j$ in~$D_j$ and~$U_\alpha \to U$ in~$\tau$.
We construct outer functions
\begin {gather*}
W_0^{(\alpha)} = \Phi_0 \left(\log \left[f_1^{\frac {1
- \alpha_1} {1 - \beta_1}} u_1^{(\alpha)}\right]\right), \quad W_1^{(\alpha)} =
\Phi_0 \left(\log \left[f_0^{\frac {\alpha_0} {\beta_0}}
u_0^{(\alpha)}\right]\right),
\\
W_0 = \Phi_0 \left(\log \left[f_1^{\frac {1
- \alpha_1} {1 - \beta_1}} u_1\right]\right), \quad W_1 = \Phi_0
\left(\log \left[f_0^{\frac {\alpha_0} {\beta_0}} u_0\right]\right)
\end {gather*}
such that~$\left|W_0^{(\alpha)}\right| = f_1^{\frac {1 - \alpha_1} {1 -
\beta_1}} u_1^{(\alpha)}$, $\left|W_1^{(\alpha)}\right| = f_0^{\frac {\alpha_0}
{\beta_0}} u_0^{(\alpha)}$, $|W_0| = f_1^{\frac {1 - \alpha_1} {1 - \beta_1}} u_1$,
and~$|W_1| = f_0^{\frac {\alpha_0} {\beta_0}} u_0$ on~$\mathbb T
\times \Omega$ using the map~$\Phi_0$ defined by~\eqref {phi1def} above.
Then~$W_j^{(\alpha)} \to W_j$ in~$\tau$, and thus also~$U_\alpha W_0^{(\alpha)}
\to U W_0$ and~$(1 - U_\alpha) W_1^{(\alpha)} \to (1 - U) W_1$ in~$\tau$.
Observe that~$\left\|U_\alpha W_0^{(\alpha)}\right\|_{F_0} \leqslant
c_2 t^{1 - \eta_0}$
and~$\left\|(1 - U_\alpha) W_1^{(\alpha)}\right\|_{F_1} \leqslant c_3
t^{-\eta_1}$, and by Proposition~\ref {latballic} we may pass to the limit in~$\tau$ in these
estimates.  This shows that
$$
(\log u_0, \log u_1, U) \in \Phi_1 (\log u_0, \log
u_1),
$$
and map~$\Phi_1$ indeed has closed graph.

Map~$\Phi_2$ is also upper semicontinuous.  This is verified in the same way
as~$\Phi_1$: if nets~$\log u_j^{(\alpha)}$, $U_\alpha$ are as above
and
$$
\left(\log v_0^{(\alpha)}, \log v_1^{(\alpha)}\right) \in \Phi_2 \left(\log
u_0^{(\alpha)}, \log u_1^{(\alpha)}, U_\alpha\right)
$$
are such that~$\log v_j^{(\alpha)} \to \log v_j$ in~$D_j$, $j \in
\{0, 1\}$, then we construct outer functions~$V_j^{(\alpha)} = \Phi_0 \left(\log
v_j^{(\alpha)}\right)$, $V_j = \Phi_0 \left(\log v_j\right)$ and pass to the
limit in the estimates
$$
\left\|U_\alpha W_0^{(\alpha)} \left[V_0^{(\alpha)}\right]^{-1}\right\|_{X_0^{1-
\zeta_0} X_1^{\zeta_0 - \gamma_1}} \leqslant 2 c_2 t^{1 - \eta_0},
$$
$$
\left\|(1 - U_\alpha) W_1^{(\alpha)}
\left[V_1^{(\alpha)}\right]^{-1} \right\|_{X_0^{1- \zeta_1 - \gamma_0}
X_1^{\zeta_1}} \leqslant 2 c_3 t^{-\eta_1}
$$ with the help of Proposition~\ref {latballic}
to show that indeed
$$
(\log v_0, \log v_1) \in \Phi_2 (\log u_0, \log u_1, U).
$$

Thus~$\Phi$ is closed and compact as a composition of compact upper
semicontinuous maps, and by Theorem~\ref {parkfptc} there exist
some~$\log u_j \in D_j$, $j \in \{0, 1\}$ such that~$(\log u_0, \log u_1) \in
\Phi (\log u_0, \log u_1)$.  This means that for some~$U \in C B_{\hclassg
{\infty}}$ we have estimates
$$
\left\|U f_1^{\frac {1 - \alpha_1} {1 - \beta_1}}\right\|_{X_0^{1-
\zeta_0} X_1^{\zeta_0 - \gamma_1}}
=
\left\|U f_1^{\frac {1 - \alpha_1} {1 - \beta_1}} u_1 u_1^{-1}\right\|_{X_0^{1-
\zeta_0} X_1^{\zeta_0 - \gamma_1}} \leqslant 2 c_2 t^{1 - \eta_0},
$$
$$
\left\|(1 - U) f_0^{\frac {\alpha_0} {\beta_0}}\right\|_{X_0^{1- \zeta_1 -
\gamma_0} X_1^{\zeta_1}} =
\left\|(1 - U) f_0^{\frac {\alpha_0} {\beta_0}}
u_0 u_0^{-1} \right\|_{X_0^{1- \zeta_1 - \gamma_0} X_1^{\zeta_1}} \leqslant 2
c_3 t^{-\eta_1}.
$$
Observe that by the choice of the parameters we have~$X_0^{1- \zeta_0}
X_1^{\zeta_0 - \gamma_1} = E_0^{\frac {1 - \alpha_1} {1 - \beta_1}}$
and~$X_0^{1- \zeta_1 - \gamma_0} X_1^{\zeta_1} = E_1^{\frac {\alpha_0}
{\beta_0}}$, and by simple computations~$\eta_0 = 1 - \frac {\beta_1 - \beta_0}
{\alpha_1 - \alpha_0} \cdot \frac {1 - \alpha_1} {1 - \beta_1}$, $\eta_1 =
\frac {\beta_1 - \beta_0} {\alpha_1 - \alpha_0} \cdot \frac {\alpha_0}
{\beta_0}$.  The estimates above imply that
$$
\left\|
|U|^{\frac {1 - \beta_1} {1 - \alpha_1}} f_1
\right\|_{E_0}
=
\left\|
U f_1^{\frac {1 - \alpha_1} {1 - \beta_1}}
\right\|_{E_0^{\frac {1 -
\alpha_1} {1 - \beta_1}}}^{\frac {1 - \beta_1} {1 - \alpha_1}}
\leqslant c_4 t^{(1 - \eta_0) \frac {1 - \beta_1} {1 - \alpha_1}} = c_4
t^{\frac {\beta_1 - \beta_0} {\alpha_1 - \alpha_0}},
$$
$$
\left\|
|1 - U|^{\frac {\beta_0} {\alpha_0}} f_0
\right\|_{E_1}
=
\left\|
(1 - U) f_0^{\frac {\alpha_0} {\beta_0}}
\right\|_{E_1^{\frac {\alpha_0} {\beta_0}}}^{\frac {\beta_0} {\alpha_0}}
\leqslant c_4 t^{-\eta_1 \frac {\beta_0} {\alpha_0}} = c_4 t^{- \frac
{\beta_1 - \beta_0} {\alpha_1 - \alpha_0}} $$
with a constant~$c_4$ independent of~$f_0$, $f_1$ and~$t$.
Since~$s = t^{\frac {\beta_1 - \beta_0} {\alpha_1 - \alpha_0}}$ is an arbitrary
positive value, by Corollary~\ref {c:tbakstabilitygen} this implies that
couple~$(E_0, E_1)$ is indeed boundedly $\AK$-stable.

Let~$0 < \delta \leqslant \frac r 2$. Couple
$$
\left(E_0^{\frac \delta 2}, E_1^{\frac \delta 2}\right) =
\left(
{\left(F^\delta\right)\strut}^{\frac 1 2} {\left[X_0^{{\delta (\beta_1 -
\beta_0)}}\right]\strut}^{\frac 1 2},
{\left(F^\delta\right)\strut}^{\frac 1 2} {\left[X_1^{\delta (\beta_1 -
\beta_0)}\right]\strut}^{\frac 1 2}
\right), \quad F = X_0^{1 - \beta_1} X_1^{\beta_0}
$$
is boundedly $\AK$-stable by Proposition~\ref {bakrscale}, which by~\cite
[Theorem~2] {kisliakov2003} implies that couple~$\left(X_0^{\delta (\beta_1 -
\beta_0)}, X_1^{\delta (\beta_1 - \beta_0)}\right)$ is $\AK$-stable.
By Corollary~\ref {akbakeq} this couple is also boundedly
$\AK$-stable, and therefore by Proposition~\ref {bakrscale} couple~$(X_0, X_1)$
is boundedly $\AK$-stable.  The proof of Theorem~\ref {t:erbakst} is complete.

A slightly different approach to the proof of Theorem~\ref {t:erbakst} is to
establish first the corresponding asymmetrical theorem with~$X_0 = Y_0$, which
leads to slightly easier computations (including the computations in Corollary~\ref
{c:tbakstabilitygen} since we only need the case~$\alpha = 1$) and then
consecutively apply it twice, first in order to extend the bounded
$\AK$-stability from the couple~$(Y_0, Y_1)$ to the couple~$(Y_0, X_1)$ under the assumptions of
Theorem~\ref {t:erbakst}, and then extend it to the entire~$(X_0, X_1)$.

\section {Proof of Theorem~\ref {goodrealinterp}}

\label {s:goodrealinterp}

First we will establish the following general version of Proposition~\ref
{akbmoreq0i} under the additional assumption that~$\Omega$ is discrete.  This,
in particular, proves~$\text{\ref {t:c:wtbmor}} \Rightarrow
\text{\ref {t:c:akst}}$.
\begin {proposition}
\label {p:bbmors}
Let~$(X, Y)$ be a couple of quasi-Banach lattices of measurable functions
on~$\mathbb T \times \Omega$ with discrete~$\Omega$ that are $r$-convex with
some~$r > 0$ satisfying the Fatou property and property~$(*)$.  Couple~$(X, Y)$
is boundedly~$\AK$-stable if and only if it
is weak-type $\BMO$-regular.
\end {proposition}
Indeed, the ``only if'' part of the proposition is done in Corollary~\ref
{akbmoreq0rb}.
Now suppose that couple~$(X, Y)$ is weak-type $\BMO$-regular.  By
Proposition~\ref {wtbmoruniv} lattices~$Y_j = \left(\lclassg {1}, (X^r)' Y^r\right)_{\theta_j, p}$ are $\BMO$-regular with some~$0 < \theta_0 <
\theta_1 < 1$, thus couple~$(Y_0, Y_1)$ is boundedly $\AK$-stable.
It satisfies the assumptions of Theorem~\ref {t:erbakst} with~$X_0 = \lclassg
{1}$ and~$X_1 = (X^r)' Y^r$, hence couple~$(X_0, X_1)$ is boundedly
$\AK$-stable.  By Proposition~\ref {bakrscale} couple~$\left(X_0^{\frac 1 2},
X_1^{\frac 1 2}\right) = \left((X^r)'^{\frac 1 2} (X^r)^{\frac 1 2}, (X^r)'^{\frac 1 2}
(Y^r)^{\frac 1 2}\right)$ is $\AK$-stable, and by \cite [Theorem~2]
{kisliakov2003} couple~$(X^r, Y^r)$ is $\AK$-stable.  By Corollary~\ref
{akbakeq} it is boundedly $\AK$-stable, and by Proposition~\ref {bakrscale}
couple~$(X, Y)$ is also boundely $\AK$-stable as claimed.

Let us show the necessity of the weak-type $\BMO$-regularity for the
$\nplus$-stability.
Suppose that under the assumptions of Theorem~\ref
{goodrealinterp} condition~\ref {t:c:npluss} is satisfied, i.~e. couple~$(X, Y)$
is $\nplus$-stable with respect to~$(\cdot, \cdot)_{\theta, q}$.
Let~$\theta_0 < \alpha_1 < \theta < \beta_1 < \theta_1$.
From the reiteration theorem it easily follows that
$\left((X,
Y)_{\alpha_1, p}, (X, Y)_{\beta_1, q}\right)$
is also $\nplus$-stable with
respect to~$(\cdot, \cdot)_{\eta_1, s}$ with~$0 < \eta_1 < 1$ satisfying~$\theta
= (1 - \eta_1) \alpha_1 + \eta_1 \beta_1$ and all~$1 \leqslant p, q \leqslant
\infty$, since
\begin {multline*}
\left[\left((X,
Y)_{\alpha_1, p}, (X, Y)_{\beta_1, q}\right)_{\eta_1, q}\right]_A = 
\left[(X, Y)_{\theta, q}\right]_A = \left(X_A, Y_A\right)_{\theta, q} =
\\ 
\left(\left(X_A,
Y_A\right)_{\alpha_1, p}, \left(X_A, Y_A\right)_{\beta_1, q}\right)_{\eta_1,
q} \subset \left(\left[(X,
Y)_{\alpha_1, p}\right]_A, \left[(X, Y)_{\beta_1, q}\right]_A\right)_{\eta_1,
q}.
\end {multline*} 
By the same reiteration theorem the lattices in this couple are also real interpolation spaces for the
couple~$\left(X^{1 - \theta_0} Y^{\theta_0}, X^{1 - \theta_1}
Y^{\theta_1}\right)$ of Banach spaces, so they are also Banach.
We may assume that~$r \leqslant 1$, and let~$0 < \delta < r$, $E =
\left(X^\delta, Y^\delta\right)_{\alpha_1, \frac p \delta}$, $F =
\left(X^\delta, Y^\delta\right)_{\beta_1, \frac q \delta}$.
By Corollary~\ref
{c:bnps} couple
$
(E, F) = \left((X, Y)_{\alpha_1, p}^\delta, (X, Y)_{\beta_1,
q}^\delta\right)
$
is also $\nplus$-stable with respect to~$(\cdot,
\cdot)_{\eta_1, \frac s \delta}$  (see also Proposition~\ref {p:riscale}).

Let~$0 < \gamma < \eta_1 < \zeta < 1$, let~$0 < \eta < 1$ be such
that~$\eta_1 = (1 - \eta) \gamma + \eta \zeta$, and let~$\alpha = (1 - \gamma)
\alpha_1 + \gamma \beta_1$, $\beta = (1 - \zeta) \alpha_1 + \zeta \beta_1$.
By~\cite [Theorem~4.7.2] {bergh} and the reiteration theorem we have
\begin {multline*}
\left[\left(\left(X^\delta, Y^\delta\right)_{\alpha, \frac s \delta},
\left(X^\delta, Y^\delta\right)_{\beta, \frac s \delta}\right)_\eta\right]_A =
\left[\left((E, F)_{\gamma, \frac s \delta}, (E, F)_{\zeta, \frac s
\delta}\right)_\eta\right]_A =
\\
\left[(E, F)_{\eta_1, \frac s \delta}\right]_A = 
\left(E_A, F_A\right)_{\eta_1, \frac s \delta} =
\left(\left(E_A, F_A\right)_{\gamma, \frac s \delta}, \left(E_A,
F_A\right)_{\zeta, \frac s \delta}\right)_\eta \subset
\\
\left(\left[(E, F)_{\gamma, \frac s \delta}\right]_A, \left[(E, F)_{\zeta, \frac
s \delta}\right]_A\right)_\eta =
\left(\left[\left(X^\delta, Y^\delta\right)_{\alpha, \frac s \delta}\right]_A,
\left[\left(X^\delta, Y^\delta\right)_{\beta, \frac s
\delta}\right]_A\right)_\eta.
\end{multline*}
That is, couple~$\left((X, Y)_{\alpha, \frac s \delta}, (X, Y)_{\beta, \frac s
\delta}\right)$
it is $\nplus$-stable with respect to the
complex interpolation~$(\cdot, \cdot)_{\theta}$ as well.
This is a couple of lattices that are~$\frac r \delta$-convex, so by~\cite
[Theorem~5.12] {kalton1994} and the remark after it this couple
is~$\BMO$-regular.  Raising it to the power~$\delta$ yields condition~\ref
{t:c:isbmor}.  This proves~$\ref {t:c:npluss} \Rightarrow \ref
{t:c:isbmor}$.

Now we establish the equivalence in Theorem~\ref {goodrealinterp} starting
with~\ref {t:c:akst}, first ignoring the ``for all'' parts of the conditions.
Transition~$\ref {t:c:akst} \Rightarrow \ref {t:c:npluss}$ is trivial, and 
transition~$\ref {t:c:npluss} \Rightarrow \ref {t:c:isbmor}$ was just verified
above.  Transitions~$\ref {t:c:isbmor} \Rightarrow \ref {t:c:xybmor} \Rightarrow
\ref {t:c:efbmor}$ are also trivial.
Transition~$\ref {t:c:efbmor} \Rightarrow \ref {t:c:efakst}$ follows from
Proposition~\ref {p:bbmors}.
Transition~$\ref {t:c:efakst} \Rightarrow \ref {t:c:isakst}$ easily follows from
the reiteration theorem (that gives $(X, Y)_{\gamma, p} = (E, F)_{\eta, p}$
for all~$\theta_0 < \gamma < \theta_1$ with some~$0 < \eta < 1$) and~\cite
[Lemma~1.1] {kisliakov1999}.
By the same reasoning, if condition~\ref {t:c:isakst} is satisfied for
some~$0 < \alpha < \theta < \beta < 1$ then it is also satisfied for
some~$\theta_0 < \alpha < \theta < \beta < \theta_1$, and by Corollary~\ref {akbakeq} the couple in
condition~\ref {t:c:isakst} is boundedly $\AK$-stable.  By Theorem~\ref
{t:erbakst} couple~$(X, Y)$ is then boundedly $\AK$-stable.  Not only does it
prove $\ref {t:c:isakst} \Rightarrow \ref {t:c:akst}$, which completes the chain
and shows that the first~$7$ conditions of Theorem~\ref {goodrealinterp} are equivalent, but it also shows
that they are equivalent to a stronger version of condition~\ref {t:c:akst}
stating that~$(X, Y)$ is boundedly $\AK$-stable, which by Proposition~\ref
{p:bbmors} implies~$\ref {t:c:akst} \Rightarrow \ref {t:c:wtbmor}$.

Thus, we have verified the equivalence of the first 8 conditions. 
Transition~$\ref {t:c:akst} \Rightarrow \ref {t:c:xyincl}$ was discussed before Proposition~\ref
{p:xyincl}.  To get the converse, let~$X_j = X^{1 - \theta_0} Y^{\theta_0}$, $j
\in \{0, 1\}$.  These lattices are Banach by the assumptions.
Then~$X^{1 - \theta} Y^\theta = X_0^{1 - \eta} X_1^\eta$ with some~$0 < \eta <
1$ satisfying~$\theta = (1 - \eta) \theta_0 + \eta \theta_1$.
On the other hand,
$(X_A, Y_A)_{\theta, \infty} = \left((X_A, Y_A)_{\theta_0, 1}, (X_A,
Y_A)_{\theta_1, 1}\right)_{\eta, \infty} \subset \left([X_0]_A,
[X_1]_A\right)_{\eta, \infty}$ by the reiteration theorem, and therefore
$\left(X_0^{1 - \eta} X_1^\eta\right)_A \subset \left([X_0]_A,
[X_1]_A\right)_{\eta, \infty}$.  Applying Proposition~\ref {p:xyincl} to
couple~$(X_0, X_1)$ yields its $\AK$-stability, which proves~$\ref {t:c:xyincl}
\Rightarrow \ref {t:c:efakst}$.

It remains to verify the ``for all'' parts of the conditions of Theorem~\ref
{goodrealinterp}.
Suppose that couple~$(X, Y)$ is $\AK$-stable and~$(E, F)$ are as in condition~\ref
{t:c:efakst}.  Let~$Y_j = (E, F)_{\eta_j, p}$, $j \in \{0, 1\}$ with
some~$0 < \eta_0 < \eta_1 < 1$ and~$p > 0$.  Then by the reiteration
theorem~$Y_j = (X, Y)_{\zeta_j, p}$, $j \in \{0, 1\}$
with~$\zeta_0 = (1 -
\eta_0) \alpha + \eta_0 \beta$ and~$\zeta_1 = (1 - \eta_1) \alpha + \eta_1
\beta$.
We choose~$\eta_j$ so that~$\theta_0 < \zeta_0 <
\theta < \zeta_1 < \theta_1$.
Then~$(Y_0, Y_1)$ is a couple of Banach lattices, it is $\AK$-stable by~\cite
[Lemma~1.1] {kisliakov1999}, and thus boundedly $\AK$-stable by Corollary~\ref
{akbakeq}.  Applying Theorem~\ref {t:erbakst} with~$X_0 = E$, $X_1 = F$ yields
the bounded $\AK$-stability of~$(E, F)$, and hence condition~\ref {t:c:efakst}. 
The
weak-type $\BMO$-regularity of~$(E, F)$, which is condition~\ref {t:c:efbmor},
follows from Proposition~\ref {p:bbmors}.  In particular, this also yields
conditions~\ref {t:c:isakst} and~\ref {t:c:xybmor}.

Finally, observe that~$\left(X^\delta, Y^\delta\right)$ is boundedly
$\AK$-stable by Proposition~\ref {bakrscale}, and it is a couple of Banach lattices for small
enough~$\delta > 0$.
We may repeat the proof of~$\ref {t:c:npluss} \Rightarrow \ref
{t:c:isbmor}$ above for this couple with~$\frac p \delta$ and~$\frac q \delta$
in place of~$p$ and~$q$, respectively, which allows us to take any~$0 < \alpha <
\beta < 1$ and yields the $\BMO$-regularity of the couple
$
\left(\left(X^\delta, Y^\delta\right)_{\alpha, \frac p \delta}, \left(X^\delta,
Y^\delta\right)_{\beta, \frac q \delta}\right).
$
Raising it to the power~$\frac 1 \delta$ by Proposition~\ref {p:riscale} yields
condition~\ref {t:c:isbmor}.

As a final remark, we explicitly note a more general version of Corollary~\ref
{akbakeq} implied by the proof of Theorem~\ref {goodrealinterp}.
\begin {corollary}
\label {c:akbakeq}
Suppose that~$X$ and~$Y$ are quasi-normed lattices of measurable functions
on~$\mathbb T \times \Omega$ with discrete~$\Omega$ satisfying the Fatou
property and property~$(*)$.  Suppose also that~$X^{1 - \theta_j} Y^{\theta_j}$,
$j \in \{0, 1\}$ are Banach lattices with some~$0 < \theta_0 <\theta_1 < 1$.
Then~$(X, Y)$ is $\AK$-stable if and only if it is boundedly~$\AK$-stable.
\end {corollary}


{}

\begin {bibliographyblock}

\normalsize
\baselineskip=17pt

\bibliographystyle {acmx}

\bibliography {bmora}

\end {bibliographyblock}

\end{document}